\documentclass[11pt,a4paper,twoside]{article}
\usepackage{epsfig}
\usepackage[latin1]{inputenc}
\usepackage{amsmath}
\usepackage{amssymb}
\usepackage{psfrag}
\def\dsp{\displaystyle}

\def\R{\mathbb{R}}

\def\eps{\varepsilon}

\def\dz{\partial_z}

\def\dZ{\partial_Z}

\def\div{\mathrm{div}}

\newcommand{\pr}{{\bf \textit{Proof: }}}

\newcommand{\cqfd}{{\nobreak\hfil\penalty50\hskip2em\hbox{}\nobreak\hfil
$\square$\qquad\parfillskip=0pt\finalhyphendemerits=0\par\medskip}}

\newtheorem{theo}{Theorem}[section]
\newtheorem{lem}{Lemma}[section]
\newtheorem{prop}{Proposition}[section]

\newtheorem{remark}{Remark}[section]

\begin{document}

%%%%%%%%%%%%%%%%%%%%%%%%%%%%%%%%%%%%%%%%%%%%%%%%%%%%%%%%%%%%%%%%%%%%%%%%%%%%%%
%%% TITLE - NAME - ABSTRACT
%%%%%%%%%%%%%%%%%%%%%%%%%%%%%%%%%%%%%%%%%%%%%%%%%%%%%%%%%%%%%%%%%%%%%%%%%%%%%%

\begin{center}
{ \Large \bf Convergence to the Reynolds approximation \\[0.1cm] with a double effect of roughness}
\end{center}

\begin{center} \bf
\large Catherine
Choquet\footnote{Faculté des Sciences et Techniques de Saint-Jérôme - Université P. Cézanne - LATP UMR 6632, 13397 Marseille Cedex 20, France\\ Email: c.choquet@univ-cezanne.fr}, Laurent
Chupin\footnote{Université de Lyon - INSA de Lyon - Pôle de Mathématiques - CNRS, UMR5208, Institut Camille Jordan -  21~av. Jean Capelle, 69621 Villeurbanne Cedex, France\\ Email: laurent.chupin@insa-lyon.fr}, Marguerite
Gisclon\footnote{LAMA - CNRS UMR5127 - Universit\'e de Savoie, 73376 Le Bourget du Lac, France\\ Email: gisclon@univ-savoie.fr}
\end{center}

\medskip



\begin{abstract}
We prove that the lubrication approximation is perturbed by a non-regular roughness of the boundary. 
We show how the flow may be accelerated using adequate rugosity profiles on the bottom.
We explicit the possible effects of some abrupt changes in the profile.
The limit system is mathematically justified through a variant of the notion of two-scale convergence.
Finally, we present some numerical results, illustrating the limit system in the three-dimensional case.
\end{abstract}

%%%%%%%%%%%%%%%%%%%%%%%%%%%%%%%%%%%%%%%%%%%%%%%%%%%%%%%%%%%%%%%%%%%%%%%%%%%%%%
\section*{Introduction}\label{introduction}
%%%%%%%%%%%%%%%%%%%%%%%%%%%%%%%%%%%%%%%%%%%%%%%%%%%%%%%%%%%%%%%%%%%%%%%%%%%%%%
We study in this paper the effect of different very small domain irregularities on a thin film flow governed
by the Stokes equations.
There exist already some references on the subject.
Roughly speaking, they may be ordered into three categories. 
The first one is devoted to the study of the roughness effects on the flow in a channel: the height of the channel, denoted $h$ and depending on the horizontal component $x$, is fixed and the effect of small scales of the boundary is studied.
In such studies, the height of the channel is written as (see~\cite{AcPiVa,JaMi2})
$$h(x) = h_2 \left( x,\frac{x}{\eps} \right) .$$
The second category  is concerned with the specific study of the thin film assumption, that is when the height of channel is assumed to be small,  of the form (see~\cite{Dri,BeChCi,Dy,El})
$$h(x) = \eps h_1(x).$$
The third category combines the mechanical point of view of the latter lubrication studies with the analysis of the roughness  effects.
Various limit models, in special regimes, are obtained depending on the ratio between the size of the rugosities and the mean height of the domain.
In~\cite{BaCh1,BaCh,Dy,El}, the ratio is assumed to be of order one, namely
$$h(x) = \eps h_2 \left( x,\frac{x}{\eps} \right) $$
and an asymptotic analysis is performed using an homogenization process.
More recently, in~\cite{BeChCi}, the authors study the case where the narrow gap is smaller than the roughness,  namely
$$h(x) = \eps\, h_2 \left( x, \frac{x}{\eps^\alpha} \right) .$$
In all previous works, the ration between the two scales, modeled by the parameter~$\alpha$ is such that $\alpha\le 1$.
Moreover, the asymptotic model obtained in the case~$\alpha\le 1$ is always the classical Reynolds equation.
In~\cite{BrChChCoGi}, the authors consider the particular case which does not enter in  the previous framework:
$$h(x) = \eps\, h_2 \left( x, \frac{x}{\eps^2} \right) .$$
They mathematically justisfy (through a variant of the notion of two-scale convergence) that an extra term modifies the standard Reynolds equation.

\smallskip

There exists now a huge number of available data for rough surfaces, due to the increasing efficiency  of the measurements (optic or laser technics, see {\it e.g.}~\cite{BhuBla91}).
On the one hand, the existence of multiple scales of rough oscillations is for instance emphasized in~\cite{HeZhu97,ZhoLeu93}.
But no analytical solution is available for complex geometries.
These data are also often too cumbersome for  numerical simulation. 
This motivates an approach using asymptotic analysis.  
On the other hand, the tribological literature does not provide any categoric answer for the choice of characteristic roughness parameters~\cite{Kou,Wiel}.
The aim of the present paper is thus to explore  the effects due to abrupt changes in the rugosity profile.
A direct application of this kind of study is the production of nano-scale electronic  components which are nowadays formed by self assembly. 
In such a process, the solvent and the monomers are confined in a thin channel.
The block copolymers create abrupt changes of the geometry (see for instance~\cite{ChaFit}).

The case of a change of speed of order one could be easily deduced from our previous work  \cite{BrChChCoGi}.
We thus focus on more abrupt changes in the rough profile.
In particular, we aim to compare the speed of the change, assumed of order $\eps^\alpha$, $\alpha >0$, and those due to the characteristic speed of the rough oscillations, assumed of order $\eps^2$, whereas~$\eps$ corresponds to the caracteristic height of the channel.
We restrict ourself to the case where $\alpha \le 2$ to ensure that our limit model respects the structure of the classical Reynolds approximation \cite{Rey86}.
%The reader may check with some straightforward computations that the assumption $\alpha \le 2$ is %essential to ensure that the pressure remains below the order $\mathcal{O}(\eps^{-2})$.

\begin{figure}[htbp]
\centering
{\psfrag{a}{$\mathcal O(\eps^2)$}\psfrag{b}{$\mathcal O(\eps^2)$}\psfrag{c}{$\mathcal O(\eps)$}\psfrag{1}{$\mathcal O(1)$}\psfrag{2}{$\mathcal O(\eps^\alpha)$}
\includegraphics[height=3.5cm]{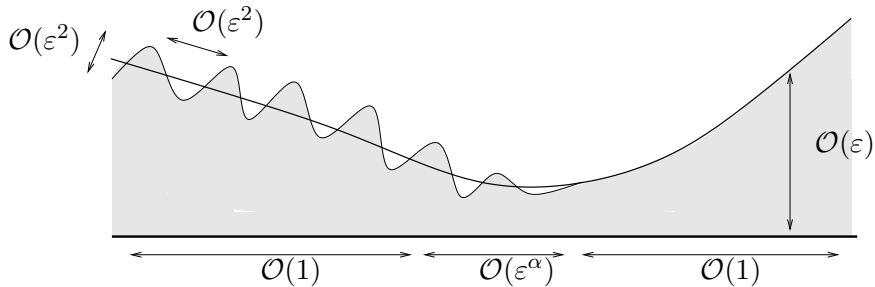}
}
\caption{Domain $\Omega^\eps$ and the different scales}\label{flow}
\end{figure}

We thus consider a model profile of the form
\begin{equation*}
h(x)=
\eps h^\eps=
\eps \left( h_1(x)+\eps \left(1-\psi\left(\frac{x}{\eps^\alpha}\right)\right) h_2\left( \frac{x}{\eps^2} \right) \right), \quad 0 \leq \alpha \le 2.
\end{equation*}
Function $\psi$ is introduced to mimic the change in the profile (see the Figure~\ref{flow}).

The domain occupied by the fluid is then:
$$\Omega^\eps = \{(x,z) \in \R^{d-1} \times \R, \quad x \in \omega, \quad 0 < z < h(x)\}$$
where $\omega$ is a domain of $\R^{d-1}$.
We choose $d=2$ or $d=3$ for the current applications.
The change in the profile occurs in a subdomain $w_T^\eps$ of~$\omega$. The domain $\omega$ is thus decomposed in $\omega= \omega_-^{\varepsilon} \cup  w_T^\eps \cup \omega_+^{\varepsilon} $.
We assume that the function $\psi$ belongs to $L^\infty(\omega)$ and is such that $\psi(x)=0$ if $ x \in  \omega_-^{\varepsilon}$ and $\psi(x)=1$ if $ x \in  \omega_+^{\varepsilon}$.
The positive function $h_1$ is the main order part of the roughness, while function $h_2$ describes the oscillating part.
For sake of simplicity, we assume that the fonction~$h_2^\eps$ defined by $h_2^\eps(x)=h_2(x/\eps^2)$ is an $\eps^2$-periodic function.
Actually, it would be sufficient for our needs to assume that $h_2$ is an admissible test function for the two-scale convergence (see Proposition \ref{prop1} below).
The first step of our method consists in introducing a new vertical variable. 
This change of variable motivates our last technical assumptions for the rough profile.
We assume that $\psi \in \mathcal C^2(\mathbb{R}^{d-1})$, $h_1 \in H^2(\mathbb{R}^{d-1})$ and $h_2^\eps \in \mathcal C^1(\mathbb{T}^{d-1})$.

\bigskip

Let us present our results.
We prove that the change in the rough profile gives the same type of correction at the main order than the oscillating part $h_2$, but only under an additional assumption linking the behavior of the roughness jumps with the one of the oscillating  part of the profile. 
Without this additional assumption, the Reynolds approximation contains a non-explicit contribution (a kind of "strange term coming from nowhere", see~\cite{CiMu}).

\begin{theo}\label{main}
Let $\omega$ be the horizontal projection of the domain of study.\\
(i) Let us denote by $\psi^0 \in L^2(\omega;H^1(\mathbb{T}^{d-1}))$ (respectively $\psi^0 \in L^2(\omega)$ if $\alpha <2$) the two-scale limit of the jump approximation $\psi^\epsilon(x)=\psi(x/\eps^\alpha)$ as $\eps \to 0$ and by $U_b$ the velocity of the lower surface.
The behavior of the thin flow is approximated by the following modified Reynolds equation on the pressure~$p$ of the fluid which only depend on the horizontal variable:
\begin{equation*}
\div_x \Bigl( \frac{h_1^3}{12 } A  \nabla_x p^0  \Bigr) = \div_x \bigl(h_1 (B U_b - Q^0) \bigr)
\quad \text{on $\omega$,}
\end{equation*}
where the functions $A$ and $B$ are defined on~$\omega$ by
 \begin{eqnarray*}
 && A = \frac{12}{N_\psi} \Bigl( e^{N_\psi/2} \int_0^1 e^{-N_\psi t^2/2}\, dt - 1 \Bigr) - \frac{12}{N_\psi} (e^{N_\psi/2}-1) \frac{ \int_0^1\int_0^s e^{N_\psi(s^2-t^2)/2}\, dtds}{ \int_0^1 e^{N_\psi s^2/2}\, ds },
 \\
 && B= \frac{1}{N_\psi}  (e^{N_\psi/2}-1) \frac{1}{ \int_0^1 e^{N_\psi s^2/2}\, ds },
 \end{eqnarray*}
where $\displaystyle N_\psi (x)= \int_{\mathbb{T}^{d-1}} \vert \nabla_X ((1-\psi^0(x,X)) h_2(X)) \vert^2 dX$, $x \in \omega$. 
The "strange" function $Q^0$ is non explicit.\\
(ii) Assuming moreover one of the following assumption,

(H1) $\nabla(\eps^2\psi^\eps h_2^\eps)$ strongly two-scale converges to $\nabla_X ( \psi^0 h_2)$,

(H2) $\psi^\eps$ strongly two-scale converges to $\psi^0$,
\\
we recover a completely explicit Reynolds-type approximation:
\begin{equation*}
\div_x \Bigl( \frac{h_1^3}{12 } A \nabla_x p^0  \Bigr) = \div_x \bigl( h_1 B U_b \bigr)\quad \text{on $\omega$.}
\end{equation*}

\end{theo}
\begin{remark}
It is crucial to clearly separate the oscillating profile $h_2$ and the jumps characteristic function $\psi$ in such a study. 
Indeed, the limit behavior described in Theorem \ref{main} clearly depends on the link between the behavior of the roughness jumps with the one of the oscillating  part of the profile (see Assumptions {\it (H1)} and {\it (H2)}). 
\end{remark}
\begin{remark}\label{base-num}
Let us describe a situation allowing the use of Assumptions {\it (H1)}.
Assume that $\omega_T^\eps$ is a finite collection of disconnected subsets of $\omega$ with measure of order $\eps^\alpha$.
As $\eps \to 0$, $\omega_T^\eps$ reduces to a finite discrete set of points $\{x_T^i\}_{i=1..N}$.
Then
$$\nabla(\eps^2\psi^\eps h_2^\eps) = \sum_{i=1}^N \bigl( \chi_{\omega_+^{i,\eps}}  \nabla_X h_2^\eps + \chi_{\omega_T^{i,\eps}}   \eps^{2-\alpha} \nabla_x \psi^\eps h_2^\eps\bigr).$$
Since the Lebesgue measure does not see accidents, one easily checks with the definition of strong two-scale convergence in Proposition \ref{prop1} that
$$\nabla(\eps^2\psi^\eps h_2^\eps) {\stackrel{2}{\to}} \sum_{i=1}^N  \chi_{\omega_+^i}  \nabla_X h_2.$$
This situation is illustrated through numerical exemples in Section \ref{SectNum}.
\end{remark}
\begin{remark}\label{rem-1202}
We recall that the usual Reynolds approximation (corresponding to $h_2=0$) reads
\begin{equation*}
\div_x \Bigl( \frac{h_1^3}{12 } \nabla_x p^0 \Bigr) = \div_x \bigl( \frac{h_1}{2} U_b \bigr)
\quad \text{on $\omega$.}
\end{equation*}
The latter result is thus a perturbation of the classical Reynolds equation.
The low order perturbations of the rough profile (oscillations $h_2$ \emph{and} abrupt change in the roughness~$\psi$) both give a perturbation of the Couette component of the Reynolds model.
This is a highly desirable effect for the applications of such a lubrication model.
We explicitly recover a perturbation which was detected with formal computations by~\cite{My}, or more recently by~\cite{PPL}.
\end{remark}
\begin{remark}
The strong two-scale convergence mentioned in {\it (ii)} of Theorem \ref{main} could appear as a rather technical assumption.
In fact, it means that the oscillation spectrum  of the sequence belongs to the integer grid (see \cite{Zhik04} for a proof through Fourier analysis).
The oscillation spectrum is then due only to the periodicity of the coefficients. 
\end{remark}  

\smallskip

\noindent 
The paper is organized as follows. In the first section, we give the model and the first estimates for the velocity and the pressure. In the second section we give the definition and the properties of a convenient two-scale convergence. 
The third section is devoted to the rigorous study of the upscaling process in view of proving Theorem \ref{main}.
In Section \ref{SectNum}, we present some numerical results, illustrating the limit system in the three-dimensional lubrication context.

%%%%%%%%%%%%%%%%%%%%%%%%%%%%%%%%%%%%%%%%%%%%%%%%%%%%%%%%%%%%%%%%%%%%%%%%%%

\section{Model and estimates}

%%%%%%%%%%%%%%%%%%%%%%%%%%%%%%%%%%%%%%%%%%%%%%%%%%%%%%%%%%%%%%%%%%%%%%%%%%

\subsection{Stokes model in lubricant context}

%%%%%%%%%%%%%%%%%%%%%%%%%%%%%%%%%%%%%%%%%%%%%%%%%%%%%%%%%%%%%%%%%%%%%%%%%%

The vector field $U^\eps=(u^\eps,w^\eps)\in \R^2\times \R$ (which describes the  fluid velocity) and the pressure (given by the scalar function $p^\eps$) satisfy the stationary Stokes equations\footnote{In all this document, the operators (like~$\Delta$) without index denote the operators with respect all the variables. To specify the operators with respect only one variable we use subscript notations, see for instance~$\Delta_x$ or~$\dz$.}:

\begin{eqnarray}
&& - \Delta u^\eps + \nabla_x p^\eps  = 0 \quad \text{in}\  \Omega^\eps,
\label{1}\\
&& - \Delta w^\eps + \dz p^\eps  = 0 \quad \text{in}\  \Omega^\eps,
\label{1b}\\
&& \div_x u^\eps + \dz w^\eps =0 \quad \text{in} \  \Omega^\eps,
\label{2}\\
&& U^\eps = U^\eps_b  \quad \text{on} \  \partial\Omega^\eps.
\label{3}
\end{eqnarray}

%\begin{eqnarray}
%&& - \Delta U^\eps + \nabla p^\eps  = 0 \quad \text{in}\  \Omega^\eps,
%\label{1}\\
%&& \div U^\eps =0 \quad \text{in} \  \Omega^\eps,
%\label{2}\\
%&& U^\eps = U^\eps_b  \quad \text{on} \  \partial\Omega^\eps.
%\label{3}
%\end{eqnarray}
\begin{itemize}
\item The viscosity of the fluid is set equal to~$1$ for sake of simplicity.
\item In classical lubrication framework, the boundary conditions are the following: 
$$U^\eps_b(x,z) = \left\{
\begin{array}{ll}
 (u_b,0) & x \in \omega, ~ z=0, \\
 (\eps\, \widetilde{u_b}(z),0) & x \in \partial \omega,\\
 (0,0) & x \in \omega, ~ z= \eps h^\eps,
\end{array} \right.
$$
where~$u_b$ is a constant corresponding to the imposed velocity at the bottom of the mechanism, and~$\widetilde{u_b}$ is a regular fonction (for instance an affine function) connecting $u_b$ for $z=0$ and $0$ for $z= h^\eps$.
Notice that for such a problem (that is for an asymptotic study $\eps\rightarrow 0$) the precise value for the velocity~$\widetilde{u_b}$ is not primary.
In fact, the physical quantity which persists to the limit is the total flux~$\int_{\partial \omega} \int_0^{\eps h^\eps} \widetilde{u_b}(z)\, dz$, see~\cite{BaCh1}.
\item Finally, the pressure is normalized by $\displaystyle \int_{\Omega^\eps} p^\eps =0$.
\end{itemize}

We introduce a new vertical variable $Z$ defined by
$Z = \displaystyle \frac{z}{h}$
to consider the same problem in the fix domain $\Omega = \omega \times (0,1)$.
Problem \eqref{1}-\eqref{3} now reads\footnote{Rigorously, the function~$u^\eps$ defined in the previous Stokes system~\eqref{1}-\eqref{3} is not egal to the function~$u^\eps$ used in this rescale domain. However, to not introduce numerous notations, we denote always the same.}
\begin{eqnarray}
&& - \Delta_x u^\eps + \Bigl( \frac{\Delta h}{h} - \frac{\vert \nabla h \vert^2}{{h}^2} \Bigr) Z \partial_Z u^\eps + 2\frac{\nabla h}{h} \cdot Z \nabla_x \partial_Z u^\eps -  \frac{\vert \nabla h \vert^2}{{h}^2} Z^2 \partial_Z^2 u^\eps  -  \frac{1}{{h}^2}  \partial_Z^2 u^\eps
\nonumber \\
&& \hspace{9cm}
 + \nabla_x p^\eps - \frac{\nabla h}{h} Z \partial_Z p^\eps  = 0 \qquad \text{in}\  \Omega,
\label{4}\\
&& - \Delta_x w^\eps + \Bigl( \frac{\Delta h}{h} - \frac{\vert \nabla h \vert^2}{{h}^2} \Bigr) Z \partial_Z w^\eps + 2\frac{\nabla h}{h} \cdot Z \nabla_x \partial_Z w^\eps -  \frac{\vert \nabla h \vert^2}{{h}^2} Z^2 \partial_Z^2 w^\eps  -  \frac{1}{{h}^2}  \partial_Z^2 w^\eps
\nonumber \\
&& \hspace{11cm}
 +\frac{1}{h}  \partial_Z p^\eps  = 0 \qquad \text{in}\  \Omega,
\label{5}\\
&& \div_x u^\eps - \frac{\nabla h}{h} \cdot Z \partial_Z u^\eps + \frac{1}{h} \partial_Z w^\eps=0 \qquad \text{in} \  \Omega,
\label{6}\\
&& U^\eps = U_b  \qquad \text{on} \  \partial\Omega,
\label{7}
\end{eqnarray}
Observe that $\displaystyle \frac{\Delta h}{h} - \frac{\vert \nabla h \vert^2}{{h}^2}=\div \left( \frac{\nabla h}{h} \right)$.

%Comme d'habitude, on suppose qu'il existe un rel\`evement $H^1$ {\it ad hoc} pour traiter la condition aux bords et obtenir les estimations...

%%%%%%%%%%%%%%%%%%%%%%%%%%%%%%%%%%%%%%%%%%%%%%%%%%%%%%%%%%%%%%%%%%%%%%%%%%
\mathversion{bold}
\subsection{Estimates and dependance with respect to the small parameter~$\eps$}\label{estimatessection}
\mathversion{normal}
%%%%%%%%%%%%%%%%%%%%%%%%%%%%%%%%%%%%%%%%%%%%%%%%%%%%%%%%%%%%%%%%%%%%%%%%%%

We begin by some estimates for the velocity.

\begin{lem}\label{lemme1}
There exists a constant $C$ such that the velocity components satisfy the following uniform estimates:
\begin{eqnarray}
&& \Vert \nabla_x U^\eps \Vert_{(L^2(\Omega))^d} \le \frac{C}{\eps},
\label{estgradu}\\
&& \Vert \partial_Z U^\eps \Vert_{(L^2(\Omega))^d} \le C,
\label{dZu} \\
&& \Vert U^\eps \Vert_{(L^2(\Omega))^d} \le C.
\label{estu}
\end{eqnarray}
\end{lem}

\noindent\pr 
These estimates are directly derived from the original problem after an adequate lifting of the non-homogenenous boundary conditions.
We write the standard energy estimate for Stokes type system and then use the change of variables to control the derivatives. We conclude with the Poincar\'e inequality.
\cqfd

We now give some uniform estimates for the pressure.

\begin{lem}
There exists a constant $C$ such that the pressure satisfies  
\begin{eqnarray}
&& \Vert \nabla_x p^\eps \Vert_{H^{-1}(\Omega)} \le \frac{C}{\eps^2},
\label{estgradp}\\
&& \Vert \partial_Z p^\eps \Vert_{H^{-1}(\Omega)} \le \frac{C}{\eps},
\label{estdZp} \\
&& \Vert p^\eps \Vert_{L^2(\Omega)} \le \frac{C}{\eps^2}.
\label{estp}
\end{eqnarray}
\end{lem}

\noindent\pr
Such estimates come from Equations (\ref{4})-(\ref{5}) estimated in  $H^{-1}$. This provides bounds on $\nabla_x p^\eps$ and $\partial_Z p^\eps$. Using Poincar\'e-Wirtinger inequality (with the normalization hypothesis) we get the $L^2$-estimates for~$p^\eps$.
\cqfd

\bigskip

\noindent {\bf Notational convention: } In what follows, for any function $\phi \in H^1(\omega \times \mathbb{T}^{d-1})$, we denote by $\nabla \phi(x,x/\eps^2) = \nabla \phi^\eps = \nabla_x \phi^\eps + \eps^{-2} \nabla_X \phi^\eps$ its horizontal gradient.

%%%%%%%%%%%%%%%%%%%%%%%%%%%%%%%%%%%%%%%%%%%%%%%%%%%%%%%%%%%%%%%%%%%%%%%%%%

\section{Definition and properties of a convenient two-scale convergence}

%%%%%%%%%%%%%%%%%%%%%%%%%%%%%%%%%%%%%%%%%%%%%%%%%%%%%%%%%%%%%%%%%%%%%%%%%%

The proof of the homogenization process will be carried out by using a variant of the
 two-scale convergence  introduced by G. Nguetseng in \cite{Ng} and
developed by G. Allaire in \cite{A}.
Let us give the basic definition and properties of this concept.

\begin{prop} \label{prop1}
A sequence  $(v^{\eps})$  of functions in $L^2(\Omega )$
two-scale converges to a limit  $v^0(x,Z,X)$ belonging to
$L^2(\Omega  \times \mathbb{T}^{d-1} )$,
$v^\eps {\stackrel{2}{\rightharpoonup}} \, v^0$,  if
$$
\lim_{\eps \to 0}\int_{\Omega}v^{\eps}(x,Z) \,
 \Psi (x,Z,x/\eps^2) \,  dx \,dZ
=\int_{\Omega }\int_{\mathbb{T}^{d-1}}v^0(x,Z,X) \,  \Psi(x,Z,X)\,  dx \,dZ\,dX,
$$
for any test function  $\Psi(x,Z,X)$, X-periodic in the third variable,
satisfying
$$
\lim_{\eps \to 0} \int_{\Omega } \vert \Psi (x,Z,x/\eps^2)
\vert^2\,  dx\,dZ
= \int_{\Omega } \int_{\mathbb{T}^{d-1}} |\Psi(x,Z,X)|^2\,  dx\,dZ\,dX .
$$
Such a function $\Psi$ is  called an admissible test function for the two-scale convergence.
Note that $Z$ is only a parameter for this definition.

\noindent (i) 
From each bounded sequence $ (v^{\eps} )$ in  $L^2 (\Omega )$
one can extract a subsequence which  two-scale converges to some limit $v^0 \in L^2(\Omega  \times \mathbb{T}^{d-1} )$. The weak $L^2$-limit of $v^\eps$ is $v(x,Z)=\int_{\mathbb{T}^{d-1}} v^0(x,Z,X) dX$.

\noindent (ii)
Let $(v^{\eps} )$ be a bounded sequence in  $L^2(0,1;H^1 (\omega))$
which converges weakly to $v$ in  $L^2(0,1;H^1 (\omega))$.
Then  $v^{\eps}{\stackrel{2}{\rightharpoonup}} \, v$  and there exists
a function $v^2 \in L^2(\Omega;H^1(\mathbb{T}^{d-1}))$ such that,
up to a subsequence,
$$\nabla v^{\eps}  {\stackrel{2}{\rightharpoonup}} \,  \nabla v (x)
+ \nabla_X v^2(x,y).$$

\noindent(iii) 
Let $ (v^{\eps} )$ be a bounded sequence in
$L^2 (\Omega)$ such that $(\eps^2 \nabla_x  v^{\eps} )$ is bounded
in $(L^2(\Omega))^{d-1}$.
 Then, there exists a function $v^0 \in L^2 (\Omega;H^1(\mathbb{T}^{d-1}))$
such that, up to a subsequence,
\begin{equation*}
v^{\eps}   {\stackrel{2}{\rightharpoonup}} \, v^0
\qquad \text{and} \qquad
\eps^2 \nabla  v^{\eps}  {\stackrel{2}{\rightharpoonup}} \nabla_X v^0 (x,Z,X ).
\end{equation*}
\end{prop}

\begin{itemize}
\item Due to our assumptions for $h_1$ and $h_2$, functions $h_1(x)$, $h_2(x/\eps^2)$, $\partial_{x_i} h_1(x)$, $\partial_{x_i}(\eps^2 h_2(x/\eps^2))=\partial_{X_i} h_2(x/\eps^2)$ can obviously be considered as admissible test functions for the two-scale convergence.
\item We also note that the function $v_2$ appearing in part {\it (ii)} of the previous proposition is the rigorous counterpart of the third term of the formal anisotropic expansion associated with the present setting:
$$ v^\eps = \tilde v^0(x,Z,\frac{x}{\eps^2}) + \eps \tilde v^1(x,Z,\frac{x}{\eps^2}) + \eps^2 \tilde v^2(x,Z,\frac{x}{\eps^2}) + \dots$$
\item An  consequence of the above definition of the two-scale convergence is the following.
\end{itemize}

\begin{lem}\label{lem1}
Let $ (v^{\eps} )$ be a bounded sequence in
$L^2 (\Omega)$ such that $(\eps^\eta \nabla  v^{\eps} )$ is bounded
in $L^2$, with $\eta<2$.
 Then  the two-scale limit $v^0 \in L^2 (\Omega ;H^1 (\mathbb{T}^{d-1} ) )$ of $v^\eps$ is such that
$$\nabla_X v^0 = 0.$$
\end{lem}

\noindent\pr
By {\it (iii)} of Proposition \ref{prop1}, 
$ \eps^2 \nabla  v^{\eps}  {\stackrel{2}{\rightharpoonup}}
\nabla_X v^0$. 
Since $(\eps^\eta \nabla  v^{\eps} )$ is bounded in $L^2$, it two-scale converges to some limit $\eta^0$ and $(\eps^2 \nabla  v^{\eps}=\eps^{2-\eta} \eps^\eta \nabla  v^{\eps})$ two-scale converges to 0.
Thus $\nabla_X v^0 = 0$.\cqfd

Note that the above lemma is a key result to treat the anisotropy of the rough profile.

%%%%%%%%%%%%%%%%%%%%%%%%%%%%%%%%%%%%%%%%%%%%%%%%%%%%%%%%%%%%%%%%%%%%%%%%%%

\section{Rigorous derivation of the limit model}

%%%%%%%%%%%%%%%%%%%%%%%%%%%%%%%%%%%%%%%%%%%%%%%%%%%%%%%%%%%%%%%%%%%%%%%%%%

\subsection{Convergence results}

%%%%%%%%%%%%%%%%%%%%%%%%%%%%%%%%%%%%%%%%%%%%%%%%%%%%%%%%%%%%%%%%%%%%%%%%%%

We infer from the  estimates derived in Subsection \ref{estimatessection} the following result.
\begin{lem}\label{lem:31}
There exist limit functions
\begin{equation*}
p^0 \in L^2 (\Omega ;L^2 (\mathbb{T}^{d-1} ) ),
\quad
u^0 \in L^2 (\Omega ;H^1 (\mathbb{T}^{d-1} ))
\quad \text{and} \quad
w^0 \in L^2 (\Omega ;H^1 (\mathbb{T}^{d-1} ))
\end{equation*}
such that
\begin{equation*}
\eps^2 p^\eps   {\stackrel{2}{\rightharpoonup}} \, p^0 ,
\quad
u^\eps {\stackrel{2}{\rightharpoonup}} \, u^0,
\quad
\nabla_X u^0 = 0,
\quad
w^\eps {\stackrel{2}{\rightharpoonup}} \, w^0
\quad \text{and} \quad
\nabla_X w^0 =0.
\end{equation*}
\end{lem}
\pr
In view of Lemma~\ref{lemme1}, we claim the existence of limits (up to a subsequence).
Moreover, using Lemma~\ref{lem1}, we assert that
$\nabla_X u^0 = 0$ and $\nabla_X w^0 =0$.
\cqfd
Before passing to the limit in the equations, we state some auxiliary results. We begin by the pressure function.
\begin{lem} \label{pressurelemma}
The two-scale limit pressure is such that
\begin{equation*}
\nabla_X p^0 = 0
\quad \text{and} \quad
\partial_Z p^0 = 0.
\end{equation*}
\end{lem}
\pr
We multiply Equation (\ref{4}) by an admissible  test function in the form $\eps^4 \phi(x,Z,x/\eps^2)$ and we integrate by parts.
Unsing in particular $\displaystyle \frac{\Delta h}{h} - \frac{\vert \nabla h \vert^2}{{h}^2}=\div \left( \frac{\nabla h}{h} \right)$, we get
\begin{eqnarray*}
&&\int_\Omega \nabla_x u^\eps \cdot ( \eps^4\nabla_x \phi^\eps + \eps^2 \nabla_X \phi^\eps) \, dx \, dZ
\\
&& \
- \int_\Omega \eps^4 (\nabla_x h_1 + \eps^{-1} (1-\psi^\eps) \nabla_Xh_2^\eps - \eps^{1-\alpha} h_2^\eps \nabla_x \psi^\eps) \cdot \frac{1}{h_1+\eps (1-\psi^\eps) h_2^\eps} \nabla_x u^\eps \partial_Z(Z\phi^\eps) \, dx \, dZ\\
&& \
- \int_\Omega \eps^4 (\nabla_x h_1 + \eps^{-1} (1-\psi^\eps) \nabla_X h_2^\eps -  \eps^{1-\alpha} h_2^\eps \nabla_x \psi^\eps) \cdot \frac{1}{h_1+\eps (1-\psi^\eps)h_2^\eps} Z \partial_Z u^\eps (\nabla_x \phi^\eps + \frac{1}{\eps^2} \nabla_X \phi^\eps) \, dx \, dZ\\
&& \
+ \int_\Omega \eps^4 \frac{\vert \nabla_x h_1 + \dsp \eps^{-1} (1-\psi^\eps) \nabla_X h_2^\eps - \dsp  \eps^{1-\alpha} h_2^\eps \nabla_x \psi^\eps \vert^2}{\vert h_1+\eps (1-\psi^\eps) h_2^\eps\vert^2} \partial_Z u^\eps \cdot  \partial_Z (Z^2 \phi^\eps) \, dx \, dZ
\\
&& \
+  \int_\Omega \eps^4 \frac{1}{\eps^2 \vert h_1+\eps (1-\psi^\eps) h_2^\eps\vert^2} \partial_Z u^\eps \cdot  \partial_Z \phi^\eps \, dx \, dZ
- \int_\Omega \eps^2 p^\eps ( \eps^2 \div_x \phi^\eps + \div_X \phi^\eps) \, dx \, dZ
 \\
&& \
+ \int_\Omega \eps^2  (\nabla_x h_1 + \eps^{-1} (1-\psi^\eps) \nabla_Xh_2^\eps -  \eps^{1-\alpha} h_2^\eps \nabla_x \psi^\eps) \cdot \frac{1}{h_1+\eps (1-\psi^\eps) h_2^\eps} \eps^2 p^\eps \partial_Z(Z\phi^\eps) \, dx \, dZ
=0.
\end{eqnarray*}
We bear in mind that the term $\eps^{1-\alpha}$ is of lower order than $\eps^{-1}$ if $0 \le \alpha <2$ and of the same order if $\alpha=2$.
Passing to the limit $\eps \to 0$, we obtain
$$\lim_{\eps \to 0} \int_\Omega \eps^2 p^\eps \div_X \phi^\eps \, dx \, dZ
= \int_\Omega \int_{\mathbb{T}^{d-1}} p^0  \div_X \phi \, dx \, dZ \, dX =0.$$
Thus $\nabla_X p^0 = 0$ and the two-scale limit and the weak $L^2$ limit of $(\eps^2 p^\eps)$ coincide. 

Now let $\phi \in L^2(\omega;H_0^1(0,1))$. 
We write
\begin{eqnarray*}
& \int_\Omega \partial_Z(\eps^2 p^\eps) \phi \, dx \, dZ &= \eps \langle  \eps\partial_Z p^\eps, \phi \rangle_{H^{-1} \times H^1} \to 0 \\
&& = - \int_\Omega \eps^2 p^\eps \partial_Z \phi \, dx \, dZ \to - \int_\Omega p^0 \partial_Z \phi \, dx \, dZ.
\end{eqnarray*}
It follows that $\partial_Z p^0 =0$.
This ends the proof of the lemma.\cqfd
We now introduce auxiliary limit functions for sake of clearness in the computations. 
We define, $\gamma^0 \in L^2(\Omega;H^1(\mathbb{T}^{d-1}))$ and $\xi^0 \in L^2(\Omega;H^1(\mathbb{T}^{d-1}))$ such that 
\begin{eqnarray*}
&& (1-\psi^\eps) u^\eps   {\stackrel{2}{\rightharpoonup}} \, \gamma^0,
\qquad \eps^2 \nabla_x ((1-\psi^\eps) u^\eps) {\stackrel{2}{\rightharpoonup}} \, \nabla_{X} \gamma^0,
 \\
&& (1-\psi^\eps)^2 u^\eps {\stackrel{2}{\rightharpoonup}} \, \xi^0, 
\qquad \eps^2 \nabla_x ((1-\psi^\eps)^2 u^\eps) {\stackrel{2}{\rightharpoonup}} \, \nabla_{X} \xi^0 .
\end{eqnarray*}
Of course, if $0 \le \alpha <2$, we have $\nabla_X \gamma^0 =0$ and $\nabla_X \xi^0=0$.
We now study the two-scale limit of the vertical velocity component.

\begin{lem} \label{verticalvelo}
The vertical velocity component is such that
$$w^\eps {\stackrel{2}{\rightharpoonup}} \, 0.$$
\end{lem}
\pr
We already mentioned that $\nabla_X w^0 =0$ (see lemma~\ref{lem:31}).
We now prove that $\partial_Z w^0 =0$.
We multiply the divergence Equation (\ref{6}) by $\eps \phi(x,Z)$ where $\phi \in H_0^1(\Omega)$.
Integrating by parts, we obtain
\begin{eqnarray}
&& - \int_\Omega \eps u^\eps \cdot \nabla_x \phi + \int_\Omega \frac{\eps}{h_1+\eps (1-\psi^\eps) h_2^\eps} \nabla_x h_1 \cdot u^\eps \partial_Z(Z\phi)
\nonumber \\
&& \qquad
+  \int_\Omega \frac{1}{h_1+\eps (1-\psi^\eps)h_2^\eps} (1-\psi^\eps) \nabla_X h_2^\eps \cdot u^\eps \partial_Z(Z\phi)
-  \int_\Omega \frac{\eps^{2-\alpha}}{h_1+\eps (1-\psi^\eps)h_2^\eps}  h^\eps_2 \nabla_x\psi^\eps \cdot u^\eps \partial_Z(Z\phi)
\nonumber \\
&& \qquad
- \int_\Omega \frac{1}{h_1+\eps (1-\psi^\eps)h_2^\eps} w^\eps \partial_Z \phi =0.
\label{dZwo2}
\end{eqnarray}
Thanks to the relation
$$ \lim_{\eps \to 0}  \int_\Omega \frac{1}{h_1+\eps (1-\psi^\eps) h_2^\eps} \nabla_X h_2^\eps \cdot (1-\psi^\eps) u^\eps \partial_Z(Z\phi)
= \int_\Omega \frac{1}{h_1} \Bigl( \int_{\mathbb{T}^{d-1}} \nabla_X  h_2 \, \gamma^0 \, dX\Bigr)  \partial_Z(Z\phi),$$
which is justified by $(1-\psi^\eps) u^\eps   {\stackrel{2}{\rightharpoonup}} \, \gamma^0$
and to the relation
\begin{eqnarray*}
&&  \!\!\!\!\!   \!\!\!\!\!  - \lim_{\eps \to 0}  \int_\Omega  \frac{\eps^{2-\alpha}}{h_1+\eps (1-\psi^\eps)h_2^\eps}  h^\eps_2 \nabla_x \psi^\eps \cdot u^\eps \partial_Z(Z\phi)
\\
&& \quad
= - \lim_{\eps \to 0}  \int_\Omega  \frac{1}{h_1+\eps (1-\psi^\eps)h_2^\eps}  h^\eps_2 \eps^2 \nabla_x (\psi(x/\eps^\alpha)  u^\eps ) \partial_Z(Z\phi)
\\
&& \qquad
+ \lim_{\eps \to 0}  \int_\Omega  \frac{1}{h_1+\eps (1-\psi^\eps)h_2^\eps}  h^\eps_2 \psi^\eps \eps^2 \div_x  u^\eps  \partial_Z(Z\phi)
\\
&& \quad
= \lim_{\eps \to 0}  \int_\Omega  \frac{1}{h_1+\eps (1-\psi^\eps)h_2^\eps}  h^\eps_2 \eps^2 \nabla_x ((1-\psi(x/\eps^\alpha) ) u^\eps ) \partial_Z(Z\phi)
\\
&& \quad
= \int_\Omega \frac{1}{h_1} \partial_Z(Z\phi) \int_{\mathbb{T}^{d-1}}  h_2 \, \nabla_X \gamma^0 \, dX, 
\end{eqnarray*}
which is justified by $\eps^2 \nabla_x ((1-\psi^\eps) u^\eps) {\stackrel{2}{\rightharpoonup}} \, \nabla_{X} \gamma^0$, 
we compute
\begin{eqnarray*}
&&   \!\!\!\!\!  \lim_{\eps \to 0}  \Bigl( \int_\Omega \frac{1}{h_1+\eps (1-\psi^\eps) h_2^\eps} \nabla_X h_2^\eps \cdot (1-\psi^\eps) u^\eps \partial_Z(Z\phi)
-  \int_\Omega  \frac{\eps^{2-\alpha}}{h_1+\eps (1-\psi^\eps)h_2^\eps}  h^\eps_2 \nabla_x \psi^\eps \cdot u^\eps \partial_Z(Z\phi)  \Bigr)
\\
&& \qquad
= \int_\Omega \frac{1}{h_1} \partial_Z(Z\phi) \int \nabla_X (h_2\gamma^0) \, dX \\
&& \qquad = 0 \quad \mbox{ because of the periodicity of } h_2\gamma^0.
\end{eqnarray*}
We thus infer from (\ref{dZwo2}) that
$$  \lim_{\eps \to 0}   \int_\Omega \frac{1}{h_1+\eps (1-\psi^\eps) h_2^\eps} w^\eps \partial_Z \phi
= \int_\Omega \frac{1}{h_1} w^0 \partial_Z \phi =0.$$
It follows that $\partial_Z w^0=0$. 
We conclude using the boundary conditions.\cqfd

\begin{remark}
We recall that the aim of this work is to get a perturbation of the Reynolds approximation.
And the previous result is characteristic of such a lubrication approximation:  it states that we can neglect the vertical component of the velocity.
%Following the lines of the previous proof, the reader may check that the assumption $\alpha>2$ does %not allow to get $w^0=0$. 
\end{remark}

%%%%%%%%%%%%%%%%%%%%%%%%%%%
%%%%%%%%%%%%%DIVERGENCE
%%%%%%%%%%%%%%%%%%%%%%%%%%%

\subsection{Divergence equation}

%%%%%%%%%%%%%%%%%%%%%%%%%%%

%
\begin{lem}\label{divlem}
We have the relation
$$h_1 \div_x  u^0  -  \nabla_x h_1  \cdot Z \, \partial_Z u^0  +  \partial_Z (w^1+L)=0,$$
where 
$w^1 = (1/d-1) \int_{\mathbb{T}^{d-1}}  \tilde w_1 \, dX$,
with $\tilde w_1\in L^2(\Omega;H^1(\mathbb{T}^{d-1}))$  is defined by the following two-scale convergence
$$ \left\{\begin{array}{l}
{\displaystyle \eps w^\eps  {\stackrel{2}{\rightharpoonup}} \, 0,}\\
{\displaystyle  \nabla (\eps w^\eps)  {\stackrel{2}{\rightharpoonup}} \,  \nabla_X \tilde w_1},
\end{array}
\right.$$
while $\partial_Z L$ is defined by the following weak convergence
$$ \left\{
\begin{array}{l}
\partial_Z L=\displaystyle \lim_{\eps \to 0}   \partial_Z \Bigl( Z \eps \displaystyle \frac{(1-\psi^\eps)h_2^\eps}{h_1} \div u^\eps \Bigr) \quad \mathrm{in}\  \mathcal{D}'(\Omega), 
\\
L_{\vert Z=0} = L_{\vert Z=1} =0.
\end{array}
\right.$$ 
\end{lem}
\begin{remark}
The term $\tilde w^1$ is the rigorous counterpart of the first order term of a formal anisotropic expansion of $w^\eps$:
$$ w^\eps = \tilde w^0(x,Z,x/\eps^2) + \eps \tilde w^1(x,Z,x/\eps^2)+ \eps^2 \tilde w^2(x,Z,x/\eps^2) + \cdots $$
\end{remark}
\pr
We multiply the divergence Equation (\ref{6}) by a test function $\phi(x,Z) \in \mathcal{D}(\Omega)$.
We obtain:
\begin{eqnarray*}
&& - \int_\Omega  u^\eps \cdot \nabla_x \phi + \int_\Omega \frac{1}{h_1+\eps (1-\psi^\eps)h_2^\eps} \nabla_x h_1 \cdot u^\eps \partial_Z(Z\phi) \, dx \, dZ
 \\
&& \quad 
+  \int_\Omega \frac{1}{h_1+\eps (1-\psi^\eps)h_2^\eps} \frac{(1-\psi^\eps)}{\eps} \nabla_X h_2^\eps \cdot u^\eps \partial_Z(Z\phi)\, dx \, dZ
 \\
&& \quad 
-  \int_\Omega \frac{1}{h_1+\eps (1-\psi^\eps)h_2^\eps}  \frac{h_2^\eps}{\eps^{\alpha-1}} \nabla_x \psi^\eps	 \cdot u^\eps \partial_Z(Z\phi)\, dx \, dZ
- \int_\Omega \frac{1}{h_1+\eps (1-\psi^\eps)h_2^\eps} \frac{1}{\eps}  w^\eps \partial_Z \phi \, dx \, dZ =0.
\end{eqnarray*}
It follows that
\begin{eqnarray}
 \int_\Omega u^0 \cdot \nabla_x \phi + \int_\Omega \frac{1}{h_1} \nabla_x h_1  \cdot Z \partial_Z u^0 \phi \, dx \, dZ  \nonumber  \\
+ \lim_{\eps \to 0} \Bigl( \int_\Omega \frac{1}{\eps} \frac{(1-\psi^\eps)\nabla_X h_2^\eps- \eps^{2-\alpha} h_2^\eps \nabla_x \psi^\eps}{h_1+\eps (1-\psi^\eps)h_2^\eps}  \cdot Z \partial_Z u^\eps \phi \, dx \, dZ
%\nonumber \\
%&& \qquad \quad
+ \int_\Omega \frac{1}{\eps} \frac{1}{h_1+\eps (1-\psi^\eps)h_2^\eps} w^\eps \partial_Z \phi  \, dx \, dZ\Bigr)  \nonumber \\ = 0.
\label{auxdiv0}
\end{eqnarray}
We note that $1/(h_1+\eps (1-\psi^\eps)h_2^\eps)$ is  a computable perturbation of $1/h_1$ for our convergences needs. 
Indeed, we have
\begin{eqnarray*}
&&  \frac{1}{h_1+\eps (1-\psi^\eps)h_2^\eps} - \frac{1}{h_1} = \frac{-\eps (1-\psi^\eps)h_2^\eps}{h_1(h_1+\eps (1-\psi^\eps)h_2^\eps)}.
\end{eqnarray*}
On the one hand, using $(1-\psi^\eps)^2 u^\eps {\stackrel{2}{\rightharpoonup}} \, \xi^0$ and $\eps^2 \nabla_x (1-\psi^\eps)^2 u^\eps {\stackrel{2}{\rightharpoonup}} \, \nabla_X \xi^0$, we write
\begin{eqnarray*}
&& \!\!\!\!\! 
\lim_{\eps \to 0} \int_\Omega \frac{1}{\eps} \frac{-\eps (1-\psi^\eps)h_2^\eps }{h_1(h_1+\eps (1-\psi^\eps)h_2^\eps)}  \bigl( (1-\psi^\eps)\nabla_X h_2^\eps - \eps^{2-\alpha}h_2^\eps \, \nabla_x \psi^\eps \bigr) \cdot Z \partial_Z u^\eps \, \phi 
\\
&& \quad
= - \int_\Omega \frac{1}{h_1^2}  \int h_2 \, \nabla_X h_2 \,  Z \, \partial_Z \xi^0 \, \phi \, dX
\\
&& \quad
- \lim_{\eps \to 0} \int_\Omega \frac{{h_2^\eps}^2 }{h_1(h_1+\eps (1-\psi^\eps)h_2^\eps)}  \frac{\eps^2}{2} \nabla_x \Bigl( (1-\psi(x/\eps^\alpha))^2 \, \partial_Z u^\eps \Bigr) \, Z \,  \phi
\\
&& \quad
+ \lim_{\eps \to 0} \int_\Omega \frac{{h_2^\eps}^2 }{h_1(h_1+\eps (1-\psi^\eps)h_2^\eps)} \frac{\eps^2}{2}(1-\psi(x/\eps^\alpha))^2 \nabla_x(  \partial_Z u^\eps ) Z \, \phi
\\
&& \quad
=- \int_\Omega \frac{1}{h_1^2}  \int h_2 \, \nabla_X h_2 \,  Z \, \partial_Z \xi^0 \, \phi \, dX
\\
&& \quad
+ \lim_{\eps \to 0} \int_\Omega \frac{{h_2^\eps}^2 }{h_1(h_1+\eps (1-\psi^\eps)h_2^\eps)}  \frac{\eps^2}{2} \nabla_x \Bigl( (1-\psi(x/\eps^\alpha))^2 \,  u^\eps \Bigr) \, \partial_Z (Z \,  \phi)
\\
&& \quad
- \lim_{\eps \to 0} \int_\Omega \frac{{h_2^\eps}^2 }{h_1(h_1+\eps (1-\psi^\eps)h_2^\eps)} \frac{\eps^2}{2}(1-\psi(x/\eps^\alpha))^2 \nabla_x(u^\eps ) \partial_Z (Z \, \phi)
\\
&& \quad
= - \int_\Omega \frac{1}{h_1^2}  \int h_2 \, \nabla_X h_2 \,  Z \, \partial_Z \xi^0 \, \phi \, dX
+ \int_\Omega \frac{1}{2h_1^2}  \int h_2^2 \, \nabla_X (\xi^0) \, \partial_Z (Z \, \phi) \, dX
\\
&& \quad
=  - \int_\Omega \frac{1}{h_1^2}  \int h_2 \, \nabla_X h_2  \, \partial_Z \xi^0 \, Z \,  \phi \, dX
- \int_\Omega \frac{1}{2h_1^2}  \int h_2^2 \, \nabla_X (\partial_Z \xi^0) \, Z \, \phi \, dX
\\
&& \quad
= - \int_\Omega \frac{1}{2h_1^2}  \int \nabla_X \bigl( h_2^2  \partial_Z \xi^0 \bigr) \,  Z \, \phi \, dX 
\\
&& \quad
=0
\end{eqnarray*}
because of the periodicity of $h_2 \xi^0$.

On the other hand, we define $w^{aux}$ by
\begin{equation}
\lim_{\eps \to 0} \int_\Omega \frac{ (1-\psi^\eps) h_2^\eps}{h_1(h_1+\eps (1-\psi^\eps)h_2^\eps)} w^\eps \partial_Z \phi 
=    \int_\Omega \Bigl( \int_{\mathbb{T}^{d-1}} w^{aux} dX \Bigr)  \partial_Z \phi =  - \int_\Omega \partial_Z \Bigl( \int_{\mathbb{T}^{d-1}}  w^{aux}  dX \Bigr)  \phi,
\label{div0}
\end{equation}
which means that
$ (1-\psi^\eps) h_2^\eps w^\eps/h_1^2 \rightharpoonup \int_{\mathbb{T}^{d-1}} w^{aux} dX$ weakly in $L^2(\Omega)$. 
Then Relation (\ref{auxdiv0}) reads:
\begin{equation}
\int_\Omega u^0 \cdot \nabla_x \phi + \int_\Omega  \frac{1}{h_1} \nabla_x h_1  \cdot Z \partial_Z u^0 \phi   + \int_\Omega \partial_Z \Bigl( \int_{\mathbb{T}^{d-1}}  w^{aux} \, dX \Bigr)  \phi =-A-B
\label{div1}
\end{equation}
where
\begin{eqnarray*}
 && A = \lim_{\eps \to 0} \int_\Omega \frac{1}{\eps} \frac{1}{h_1} \bigl( (1-\psi^\eps) \nabla_X h_2^\eps - \eps^{2-\alpha} h^\eps_2 \nabla_x \psi^\eps \bigr) \cdot Z \partial_Z u^\eps \phi,
 \\
&& B =  \lim_{\eps \to 0}  \int_\Omega \frac{1}{\eps} \frac{1}{h_1} w^\eps \partial_Z \phi.
\end{eqnarray*}
We write
\begin{eqnarray}
A = \lim_{\eps \to 0}  \int_\Omega \frac{1}{h_1} \nabla_x (\eps (1-\psi^\eps)h_2^\eps) \cdot (Z \partial_Z u^\eps) \phi
\nonumber 
\\
=  \lim_{\eps \to 0} \Bigl( 
\int_\Omega \frac{\eps (1-\psi^\eps) h_2^\eps}{h_1} \div_x u^\eps \partial_Z(Z\phi)
- \int_\Omega \frac{\eps (1-\psi^\eps)h_2^\eps}{h_1^2} \nabla_x h_1 \cdot u^\eps \partial_Z(Z\phi)
\nonumber \\
+\int_\Omega \frac{\eps (1-\psi^\eps) h_2^\eps}{h_1} u^\eps \cdot \nabla_x(\partial_Z(Z\phi))
 \Bigr)  
 \nonumber
  \\
=  \lim_{\eps \to 0}  \int_\Omega \frac{\eps (1-\psi^\eps)h_2^\eps}{h_1} \div_x u^\eps \partial_Z(Z\phi)
\nonumber  \\
=   \lim_{\eps \to 0}  \int_\Omega \frac{\eps (1-\psi^\eps)h_2^\eps}{h_1} \div_x u^\eps Z \partial_Z \phi
  + \lim_{\eps \to 0}  \int_\Omega \frac{\eps (1-\psi^\eps)h_2^\eps}{h_1} \div_x u^\eps \phi
 \nonumber\\
=  - \lim_{\eps \to 0}  \int_\Omega \partial_Z \Bigl( Z \frac{\eps (1-\psi^\eps)h_2^\eps}{h_1} \div_x u^\eps \Bigr)  \phi
  + \lim_{\eps \to 0}  \int_\Omega \frac{\eps (1-\psi^\eps)h_2^\eps}{h_1} \div_x u^\eps \phi.
\label{div2}
\end{eqnarray}
The first term in the right hand side will give the conservative contribution $\partial_Z L$ to the final divergence equation.
Let us compute the last term.
We multiply the divergence Equation (\ref{6}) by the test function $\eps (1-\psi^\eps)h_2^\eps\phi(x,Z)$. We get
\begin{eqnarray*}
&& \!\!\!\!\!\!
 \lim_{\eps \to 0} \int_\Omega \frac{ (1-\psi^\eps)h_2^\eps}{h_1+\eps (1-\psi^\eps)h_2^\eps} w^\eps \partial_Z \phi 
 \\
 && \!\!\!
=  \lim_{\eps \to 0} \int_\Omega \eps \div_x u^\eps (1-\psi^\eps)h_2^\eps \phi
+  \lim_{\eps \to 0} \int_\Omega \frac{\eps ( 1-\psi^\eps)h_2^\eps}{h_1+\eps (1-\psi^\eps)h_2^\eps} \nabla_x h_1 \cdot u^\eps \partial_Z (Z\phi)
\\
&& \  
+  \lim_{\eps \to 0} \int_\Omega \frac{ h_2^\eps}{h_1+\eps (1-\psi^\eps)h_2^\eps} \nabla_X h_2^\eps \cdot (1-\psi^\eps)^2 u^\eps \partial_Z (Z\phi)
\\
&& \  
-   \lim_{\eps \to 0} \int_\Omega \frac{ \eps^{2-\alpha}{h_2^\eps}^2}{h_1+\eps (1-\psi^\eps)h_2^\eps} \nabla_x \psi^\eps \cdot (1-\psi^\eps) u^\eps \partial_Z (Z\phi)
\\
&& \!\!\!
=  \lim_{\eps \to 0} \int_\Omega \eps \, \div_x u^\eps \, (1-\psi^\eps) \, h_2^\eps \, \phi 
+ \int_\Omega \frac{1}{2h_1} \, \int \nabla_X( h_2^2) \cdot \xi^0 \, \partial_Z(Z\phi) \, dX
\\
&& \  
+  \lim_{\eps \to 0} \int_\Omega \frac{ {h_2^\eps}^2}{h_1+\eps (1-\psi^\eps)h_2^\eps} \, \eps^2 \, \nabla_x \Bigl(\frac{1}{2} (1-\psi^\eps)^2 u^\eps \Bigr) \, \partial_Z (Z\phi)
\\
&& \  
- \lim_{\eps \to 0} \int_\Omega \frac{ {h_2^\eps}^2}{h_1+\eps (1-\psi^\eps)h_2^\eps} \frac{\eps^2}{2} (1-\psi^\eps)^2 \div_x u^\eps  \, \partial_Z (Z\phi)
\\
&& \!\!\!
=  \lim_{\eps \to 0} \int_\Omega \eps \, \div_x u^\eps \, (1-\psi^\eps) \, h_2^\eps \, \phi + \int_\Omega \frac{1}{2h_1} \int \nabla_X (h_2^2) \cdot \xi^0 \,  \partial_Z(Z\phi) \, dX + \int_\Omega \frac{1}{2h_1} \int h_2^2 \, \nabla_X \xi^0 \, \partial_Z(Z\phi) \, dX
\\
&& \!\!\! 
=  \lim_{\eps \to 0} \int_\Omega \eps \, \div_x u^\eps \, (1-\psi^\eps) \, h_2^\eps \, \phi 
+ \frac{1}{2} \int_\Omega \frac{1}{h_1} \int \nabla_X ( h_2^2  \xi^0)\,  \partial_Z(Z\phi) \, dX
\\
&& \!\!\!
=  \lim_{\eps \to 0} \int_\Omega \eps \, \div_x u^\eps \, (1-\psi^\eps) \, h_2^\eps \, \phi
\end{eqnarray*}
because of the periodicity of $h_2 \xi^0$.

Bearing in mind (\ref{div0}), we infer from the latter relation that:
\begin{gather}
  \lim_{\eps \to 0} \int_\Omega \frac{ (1-\psi^\eps) h_2^\eps}{h_1(h_1+\eps (1-\psi^\eps)h_2^\eps)} w^\eps \partial_Z \phi 
=  \lim_{\eps \to 0} \int_\Omega \frac{h_2^\eps}{h_1} \eps \, (1-\psi^\eps) \, \phi \, \div_x u^\eps 
\notag \\
= - \int_\Omega \partial_Z \Bigl( \int_{\mathbb{T}^{d-1}}  w^{aux} \, dX \Bigr) \phi
\label{div3}
\end{gather}
because of the definition of $w^{aux}$.

Let us now compute $B=\dsp \lim_{\eps \to 0}  \int_\Omega \frac{1}{\eps} \frac{1}{h_1} w^\eps \, \partial_Z \, \phi$.
We write
\begin{equation*}
\begin{aligned}
\int_\Omega \frac{1}{\eps} \frac{1}{h_1} w^\eps \partial_Z \phi 
& = \int_\Omega \frac{1}{\eps^2} (\eps w^\eps \frac{1}{h_1} \partial_Z \phi ) \\
& = \frac{1}{d-1} \int_\Omega \div (I(x/\eps^2)) (\eps w^\eps \frac{1}{h_1} \partial_Z \phi ) \\
& = - \frac{1}{d-1} \int_\Omega  I(x/\eps^2)\cdot \nabla (\eps w^\eps \frac{1}{h_1} \partial_Z \phi )
\end{aligned}
\end{equation*}
where $I$ is the $1$-periodic function such that $I(x)=x$ if $x \in \mathbb T^{d-1}$.
Then 
\begin{equation}\label{div4}
\begin{aligned}
B & =  - \frac{1}{d-1} \int_\Omega \Bigl( \int_{\mathbb{T}^{d-1}}  X \cdot \nabla_X \tilde w_1 \, dX \Bigr) \frac{1}{h_1} \,  \partial_Z \phi \\
& = \frac{1}{d-1} \int_\Omega \Bigl( \int_{\mathbb{T}^{d-1}}  \tilde w_1 \, dX \Bigr)  \frac{1}{h_1} \, \partial_Z \phi \\
& = \int_\Omega \frac{w^1}{h_1} \, \partial_Z \phi
  = - \int_\Omega \frac{\partial_Z w^1}{h_1} \, \phi
\end{aligned}
\end{equation}
where 
$w^1 = \frac{1}{d-1} \int_{\mathbb{T}^{d-1}}  \tilde w_1 \, dX$,
the function $\tilde w_1\in L^2(\Omega;H^1(\mathbb{T}^{d-1}))$ being defined by 
$$ \left\{\begin{array}{l}
{\displaystyle \eps w^\eps  {\stackrel{2}{\rightharpoonup}} \, 0,}\\
{\displaystyle  \nabla (\eps w^\eps)  {\stackrel{2}{\rightharpoonup}} \,  \nabla_X \tilde w_1.}
\end{array}
\right.$$
Using (\ref{div1})-(\ref{div4}), we obtain the announced result.
\cqfd

The conservative form of the limit divergence equation obtained in the previous lemma is sufficient to get an explicit Reynolds approximation of the flow, even if function $L$ is not explicitly computed (see \cite{BrChChCoGi}).
By the way, we end this section by computing $L$ under an additional assumption. 
\begin{lem}
Assuming moreover one of the following assumption, 

\noindent
(H1) $\nabla(\psi^\eps h_2^\eps)$ strongly two-scale converges to $\nabla_X ( \psi^0 h_2)$,\\
(H2) $\psi^\eps$ strongly two-scale converges to $\psi^0$,\\
(H3) $\alpha <2$,\\
 then $L=0$ and the limit divergence equation reads 
$$h_1 \div_x u^0 - \nabla_x h_1 \cdot Z \partial_Z u^0 
+ \partial_Z  w^1 = 0.$$
\end{lem}
\pr
From the assumptions (H1) or (H2), it follows that
$$ h_2(X) \gamma^0(x,X) = h_2 (X) \bigl( 1- \psi^0(x,X)\bigr) u_0(x,Z),$$
$$ h_2(X) \xi^0(x,X) = h_2 (X) \bigl( 1- \psi^0(x,X)\bigr)^2 u_0(x,Z).$$
The two scale limit of $Z \, \eps (1-\psi^\eps) \, h_2^\eps \, \div u^\eps$ is   $h_1L=Z \, (1-\psi^0) \, h_2 \, \div_X u^1$ where $u^1 \in L^2(\Omega;H^1(\mathbb{T}^{d-1}))$ is the anisotropic two-scale limit defined by
$$\nabla (\eps u^\eps)  {\stackrel{2}{\rightharpoonup}} \, \nabla_Xu^1.$$
Let us compute $\div_X u^1$.
To this aim, we multiply the divergence Equation (\ref{6}) by $\eps \phi(x,Z,x/\eps^2)$ where $\phi$ is an admissible test function for the two scale convergence.
Integrating by parts, we obtain
\begin{eqnarray*}
&& \int_\Omega \eps \div u^\eps \phi^\eps 
- \int_\Omega \bigl( (1-\psi^\eps)\nabla_X h_2^\eps -  \eps^{2-\alpha}h_2^\eps \nabla_x \psi^\eps \bigr) \cdot Z\partial_Z u^\eps  \frac{\phi^\eps}{h_1+\eps (1-\psi^\eps)h_2^\eps} 
\\
&& \qquad
- \int_\Omega \eps \nabla_x h_1 \cdot Z\partial_Z u^\eps  \frac{\phi^\eps}{h_1+\eps (1-\psi^\eps) h_2^\eps} + \int_\Omega \partial_Z w^\eps \frac{\phi^\eps}{h_1+\eps (1-\psi^\eps) h_2^\eps} 
=0.
\end{eqnarray*}
We recall that $\eps \partial_Z u^\eps {\stackrel{2}{\to}} \, 0$ (a strong two-scale convergence due to the fact that $\partial_Z u^\eps$ is bounded in $(L^2(\Omega))^d$) and $w^\eps {\stackrel{2}{\rightharpoonup}} \, 0$ (see Lemma \ref{verticalvelo}).
We also have
\begin{eqnarray*}
&& \lim_{\eps \to 0}  \int_\Omega \bigl( (1-\psi^\eps)\nabla_X h_2^\eps 
-  \eps^{2-\alpha} h_2^\eps \nabla_x \psi^\eps \bigr) \cdot Z\partial_Z u^\eps  \frac{\phi^\eps}{h_1+\eps (1-\psi^\eps)h_2^\eps} 
= \int_\Omega \int_{\mathbb{T}^{d-1}}  \frac{1}{h_1} \, \div_X (h_2 Z \partial_Z \gamma^0) \, \phi \, dX,
\end{eqnarray*}
because 
\begin{gather*}
\bigl( (1-\psi^\eps)\nabla_X h_2^\eps - \eps^{2-\alpha} h_2^\eps \nabla_x \psi^\eps \bigr) u^\eps
= \eps^2 \div_x \bigl( (1-\psi^\eps)h_2^\eps u^\eps \bigr) - \eps^2 (1-\psi^\eps)h_2^\eps \div_x u^\eps 
\\
{\stackrel{2}{\rightharpoonup}} \, \div_X \bigl( h_2 \gamma^0 \bigr).
\end{gather*}
We thus conclude that
\begin{equation}
\div_X u^1 = \frac{1}{h_1} \div_X (h_2 \cdot Z \partial_Z \gamma^0).
\label{divu12}
\end{equation}

With the additional assumption (H1) or (H2), we infer from (\ref{divu12}) and from $\nabla_x u^0=0$ that
$$ \div_X u^1 =  \frac{1}{h_1} \div_X (h_2 (1-\psi^0) Z \partial_Z u^0) =  \frac{1}{h_1} \nabla_X (h_2 (1-\psi^0) )\cdot Z \partial_Z u^0.$$
We conclude by computing
$$ \int_{\mathbb{T}^{d-1}}  L \, dX
=  \frac{1}{h_1^2}  \partial_Z ( Z^2 \partial_Z u_0 )\, \Bigl( \int_{\mathbb{T}^{d-1}}  h_2 \, (1-\psi^0) \, \nabla_X (h_2 (1-\psi^0))  \,  dX\Bigr) = 0.$$

Assume now (H3).
Following the lines of the proof of Lemma \ref{verticalvelo}, one easily check that the assumption $\alpha <2$ leads to
$$ (1-\psi^\eps) w^\eps {\stackrel{2}{\rightharpoonup}} \, 0.$$
Now, using the test function $\eps (1-\psi^\eps) \phi(x,Z,x/\eps^2)$ instead of $\eps \phi(x,Z,x/\eps^2)$ in the latter derivation of $\div_X u^1$, we state that
$$ \eps (1-\psi^\eps) \div u^\eps {\stackrel{2}{\rightharpoonup}} \, \frac{1}{h_1} \div_X \bigl( h_2  \cdot Z \partial_Z \xi^0 \bigr).$$
Integrating over $\mathbb{T}^{d-1}$, we recover once again $L=0$.
\cqfd

%%%%%%%%%%%%%%%%%%%%%%%%%%%
%%%%%%%%%%%%%MOMENTUM
%%%%%%%%%%%%%%%%%%%%%%%%%%%

\subsection{Momentum equation}

%%%%%%%%%%%%%%%%%%%%%%%%%%%

%
\begin{lem}\label{momentlem}
Let $q^0$ be the "strange" function defined by the following weak convergence in $L^2(\omega;H^{-2}(0,1))$:
$$ \bigl\vert \nabla \bigl( \eps^2 (1-\psi^\eps) h_2^\eps \bigr) \bigr\vert^2 Z \partial_Z (u^\eps-u^0)
+ Z \partial_Z \Bigl( Z \partial_Z \bigl( \frac{1}{2} \eps^2 p^\eps  \bigr) \Bigr) \nabla ( \eps^2 \vert (1-\psi^\eps)h_2^\eps\vert^2 ) 
\rightharpoonup  q^0. $$ 
The limit of the first momentum equation is 
$$ N_\psi Z \partial_Z u^0 - \partial_{ZZ}^2 u^0 + h_1^2 \nabla_x p^0 + q^0 =0 , $$
where
$$ N_\psi (x)= \int_{\mathbb{T}^{d-1}} \vert \nabla_X \bigl( (1-\psi^0(x,X)) h_2(X) \bigr) \vert^2 dX, \  x \in \omega.$$
\end{lem}

\begin{remark}
Assuming moreover {\it (H1)} or {\it (H2)},  that is $\nabla(\psi^\eps h_2^\eps)$ strongly two-scale converges to $\nabla_X ( \psi^0 h_2)$
or
$\psi^\eps$ strongly two-scale converges to $\psi^0$,  and using  $\partial_Z p^0=0$,  one checks that  $q^0=0$.
\end{remark}

\noindent\pr
Let us pass to the limit $\eps \to 0$ in the momentum equation.
We multiply Equation (\ref{4}) by $\eps^2 \phi(x,Z)$ where $\phi \in H_0^1(\omega; H_0^2(0,1))$ and we obtain:
\begin{eqnarray*}
&& \int_\Omega \eps^2 \nabla_x u^\eps \cdot \nabla_x \phi 
- \int_\Omega \frac{\eps}{h_1+\eps(1-\psi^\eps)h_2^\eps} (\nabla_x h_1 + 
\eps^{-1} (1-\psi^\eps) \nabla_X h_2^\eps - \eps^{1-\alpha}  h_2^\eps \nabla_x \psi^\eps) \cdot \eps \nabla_x  u^\eps \partial_Z(Z\phi)
\nonumber  \\
&& \quad
- \int_\Omega \frac{\eps^2}{h_1+\eps (1-\psi^\eps)h_2^\eps} (\nabla_x h_1 + \eps^{-1} (1-\psi^\eps) \nabla_X h_2^\eps - \eps^{1-\alpha}  h_2^\eps \nabla_x \psi^\eps) Z \partial_Z u^\eps \nabla_x \phi
\nonumber  \\
&& \quad
+ \int_\Omega \eps^2 \frac{\vert\nabla_x h_1 + \eps^{-1} (1-\psi^\eps) \nabla_X h_2^\eps - \eps^{1-\alpha}  h_2^\eps \nabla_x \psi^\eps \vert^2}{(h_1+\eps (1-\psi^\eps) h_2^\eps)^2}   \partial_Z u^\eps \cdot \partial_Z(Z ^2 \phi)
\nonumber  \\
&& \quad
+ \int_\Omega \frac{1}{(h_1+\eps (1-\psi^\eps) h_2^\eps)^2}  \partial_Z u^\eps \cdot \partial_Z \phi
- \int_\Omega \eps^2 p^\eps \div_x \phi
\nonumber   \\
&& \quad
+ \int_\Omega \frac{1}{h_1+\eps (1-\psi^\eps) h_2^\eps} \eps^2 p^\eps \bigl( \nabla_x h_1 + \eps^{-1} (1-\psi^\eps) \nabla_X h_2^\eps - \eps^{1-\alpha} h_2^\eps \nabla_x \psi^\eps \bigr) \cdot \partial_Z (Z \phi) = 0.
\end{eqnarray*}
Passing to the limit $\eps \to 0$, we get
\begin{eqnarray}
&& 
- \lim_{\eps \to 0} \int_\Omega \eps^2 \nabla \Bigl( \frac{(1-\psi^\eps)h_2^\eps}{h_1}  \partial_Z(Z\phi) \Bigr) \cdot \eps \nabla_x u^\eps
\nonumber    \\
&& \quad
+ \lim_{\eps \to 0} \int_\Omega  \frac{\vert (1-\psi^\eps) \nabla_X h_2^\eps - \eps^{2-\alpha}  h_2^\eps \nabla_x \psi^\eps \vert^2}{h_1^2}   \partial_Z u^\eps \cdot \partial_Z(Z ^2 \phi)
\nonumber    \\
&& \quad
+ \int_\Omega \frac{1}{h_1^2}  \partial_Z u^0 \cdot \partial_Z \phi
- \int_\Omega  p^0 \div_x \phi
+\lim_{\eps \to 0} \int_\Omega \eps^3 p^\eps \div \bigl( \frac{(1-\psi^\eps)h_2^\eps}{h_1} \partial_Z (Z \phi) \bigr) = 0.
\label{corr3}
\end{eqnarray}
We express the first term of the left hand side of (\ref{corr3}) using once again (\ref{4}), but for the test function $\eps^3 \dsp \frac{(1-\psi^\eps)h_2^\eps}{h_1}  \partial_Z(Z\phi)$.
For this choice of test function, we obtain:
\begin{eqnarray*}
&& \int_\Omega \nabla_x u^\eps \cdot \eps^3 \nabla \Bigl( \frac{(1-\psi^\eps)h_2^\eps}{h_1}  \partial_Z(Z\phi) \Bigr)
 \\
&& ~
= \int_\Omega \frac{\eps^3}{h_1+ \eps(1-\psi^\eps)h_2^\eps}  \bigl(\nabla_x h_1 + \eps^{-1} ((1-\psi^\eps)\nabla_X h_2^\eps - \eps^{2-\alpha} h_2^\eps \nabla_x \psi^\eps) \bigr) \cdot Z \partial_Z u^\eps \nabla \Bigl( \frac{(1-\psi^\eps)h_2^\eps}{h_1}  \partial_Z(Z\phi) \Bigr)
 \\
&& \quad
+ \int_\Omega  \frac{\eps^2}{h_1+ \eps(1-\psi^\eps)h_2^\eps} \bigl(\nabla_x h_1 + \eps^{-1}  (1-\psi^\eps)\nabla_X h_2^\eps - \eps^{2-\alpha} h_2^\eps \nabla_x \psi^\eps \bigr) \cdot \eps \nabla_x u^\eps \, 
 \partial_Z \Bigl( Z \frac{(1-\psi^\eps)h_2^\eps}{h_1}  \partial_Z(Z\phi) \Bigr)
 \\
&& \quad
- \int_\Omega \eps^2 \frac{\vert\nabla_x h_1 + \eps^{-1}  (1-\psi^\eps)\nabla_X h_2^\eps  - \eps^{1-\alpha} h_2^\eps \nabla_x \psi^\eps \vert^2}{(h_1+\eps (1-\psi^\eps) h_2^\eps)^2} \eps \partial_Z u^\eps \cdot \partial_Z  \Bigl( Z^2 \frac{(1-\psi^\eps)h_2^\eps}{h_1}  \partial_Z(Z\phi) \Bigr)
\\
&& \quad
- \int_\Omega \frac{1}{\eps^2(h_1+\eps (1-\psi^\eps) h_2^\eps)^2} \eps^3 \partial_Z u^\eps \cdot \partial_Z   \Bigl( \frac{(1-\psi^\eps)h_2^\eps}{h_1}  \partial_Z(Z\phi) \Bigr)
+ \int_\Omega \eps^3 p^\eps \div  \Bigl( \frac{(1-\psi^\eps)h_2^\eps}{h_1}  \partial_Z(Z\phi) \Bigr)
\\
&& \quad
- \int_\Omega \frac{1}{h_1+\eps (1-\psi^\eps)h_2^\eps} \eps^2 p^\eps \eps \bigl( \nabla_x h_1 + \eps^{-1} (1-\psi^\eps)\nabla_X h_2^\eps - \eps^{1-\alpha} h_2^\eps \nabla_x \psi^\eps \bigr) \cdot \partial_Z  \Bigl( Z \frac{(1-\psi^\eps)h_2^\eps}{h_1}  \partial_Z(Z\phi) \Bigr) .
\end{eqnarray*}
Passing to the limit $\eps \to 0$ in the latter relation, we get:
\begin{eqnarray}
&& \lim_{\eps \to 0} \int_\Omega \nabla_x u^\eps \cdot \eps^3 \nabla \Bigl( \frac{(1-\psi^\eps)h_2^\eps}{h_1}  \partial_Z(Z\phi) \Bigr)
\nonumber \\
&& \quad
= \lim_{\eps \to 0} \int_\Omega \frac{\eps^2}{h_1}  \bigl(  (1-\psi^\eps)\nabla_X h_2^\eps - \eps^{2-\alpha} h_2^\eps \nabla_x \psi^\eps \bigr) \cdot Z \partial_Z u^\eps \nabla \Bigl( \frac{(1-\psi^\eps)h_2^\eps}{h_1}  \partial_Z(Z\phi) \Bigr)
\nonumber \\
&& \quad
+ \lim_{\eps \to 0} \int_\Omega \eps^3 p^\eps \div  \Bigl( \frac{(1-\psi^\eps)h_2^\eps}{h_1}  \partial_Z(Z\phi) \Bigr)
\nonumber \\
&& \quad
- \lim_{\eps \to 0} \int_\Omega \frac{1}{h_1} \eps^2 p^\eps  \bigl(   (1-\psi^\eps)\nabla_X h_2^\eps - \eps^{2-\alpha} h_2^\eps \nabla_x \psi^\eps \bigr) \cdot \partial_Z  \Bigl( Z \frac{(1-\psi^\eps)h_2^\eps}{h_1}  \partial_Z(Z\phi) \Bigr) .
\label{corr2}
\end{eqnarray}
Inserting the previous computation  in (\ref{corr3}), we obtain:
\begin{eqnarray*}
&& 
- \lim_{\eps \to 0} \int_\Omega \frac{1}{h_1} \bigl( (1-\psi^\eps)\nabla_X h_2^\eps - \eps^{2-\alpha} h_2^\eps \nabla \psi^\eps \bigr) \cdot Z \partial_Z u^\eps \nabla \bigl(\eps^2 \frac{(1-\psi^\eps)h_2^\eps}{h_1} \partial_Z (Z \phi) \bigr)
\\
&& \quad
- \lim_{\eps \to 0} \int_\Omega \eps^3 p^\eps \div \bigl( \frac{(1-\psi^\eps)h_2^\eps}{h_1} \partial_Z (Z \phi) \bigr)
\\
&& \quad
+ \lim_{\eps \to 0} \int_\Omega \frac{1}{h_1} \eps^2 p^\eps \bigl( (1-\psi^\eps)\nabla_X h_2^\eps - \eps^{2-\alpha} h_2^\eps \nabla \psi^\eps \bigr)  \frac{(1-\psi^\eps)h_2^\eps}{h_1} \partial_Z (Z \partial_Z (Z\phi)) 
\\
&& \quad
+ \lim_{\eps \to 0} \int_\Omega  \frac{\vert (1-\psi^\eps) \nabla_X h_2^\eps - \eps^{2-\alpha}  h_2^\eps \nabla \psi^\eps \vert^2}{h_1^2}   \partial_Z u^\eps \cdot \partial_Z(Z ^2 \phi)
 \\
&& \quad
+ \int_\Omega \frac{1}{h_1^2}  \partial_Z u^0 \cdot \partial_Z \phi
- \int_\Omega  p^0 \div_x \phi
   \\
&& \quad
+\lim_{\eps \to 0} \int_\Omega \eps^3 p^\eps \div \bigl( \frac{(1-\psi^\eps)h_2^\eps}{h_1} \partial_Z (Z \phi) \bigr)=0.
\end{eqnarray*}
We thus note that the combined use of (\ref{corr3}) and (\ref{corr2}), that is the use of the momentum equation (\ref{4}) for two different {\it ad hoc} test functions, is sufficient to pass to the limit without studying the limit behavior of the non-bounded sequence $(\eps p^\eps)$. 
This is a subsequent improvement of \cite{BrChChCoGi}.
At this step, we have proven that:
\begin{eqnarray*}
&& 
- \lim_{\eps \to 0} \int_\Omega \frac{1}{h_1^2} \bigl\vert (1-\psi^\eps)\nabla_X h_2^\eps - \eps^{2-\alpha} h_2^\eps \nabla \psi^\eps \bigr\vert^2  Z \partial_Z u^\eps   \partial_Z (Z \phi) 
\\
&& \quad
+ \lim_{\eps \to 0} \int_\Omega \frac{1}{h_1} \eps^2 p^\eps \bigl( (1-\psi^\eps)\nabla_X h_2^\eps - \eps^{2-\alpha} h_2^\eps \nabla \psi^\eps \bigr)  \frac{(1-\psi^\eps)h_2^\eps}{h_1} \partial_Z (Z \partial_Z (Z \phi)) 
\\
&& \quad
+ \lim_{\eps \to 0} \int_\Omega  \frac{\vert (1-\psi^\eps) \nabla_X h_2^\eps -  \eps^{2-\alpha} h_2^\eps \nabla \psi^\eps \vert^2}{h_1^2}   \partial_Z u^\eps \cdot \partial_Z(Z ^2 \phi)
 \\
&& \quad
+ \int_\Omega \frac{1}{h_1^2}  \partial_Z u^0 \cdot \partial_Z \phi
- \int_\Omega  p^0 \div_x \phi
 =0,
\end{eqnarray*}
that is, combining the first and third lines of the latter relation,
\begin{eqnarray}
&& 
 \lim_{\eps \to 0} \int_\Omega \frac{1}{h_1} \eps^2 p^\eps \bigl( (1-\psi^\eps)\nabla_X h_2^\eps - \eps^{2-\alpha} h_2^\eps \nabla \psi^\eps \bigr)  \frac{(1-\psi^\eps)h_2^\eps}{h_1} \partial_Z (Z \partial_Z (Z \phi)) 
\nonumber \\
&& \quad
+ \lim_{\eps \to 0} \int_\Omega  \frac{\vert (1-\psi^\eps) \nabla_X h_2^\eps -  \eps^{2-\alpha} h_2^\eps \nabla \psi^\eps \vert^2}{h_1^2}  Z  \partial_Z u^\eps \cdot  \phi
 \nonumber \\
&& \quad
+ \int_\Omega \frac{1}{h_1^2}  \partial_Z u^0 \cdot \partial_Z \phi
- \int_\Omega  p^0 \div_x \phi
 =0.
 \label{presque}
\end{eqnarray}
It remains to exhibit the partial differential equation corresponding to the latter weak formulation.
Lemma \ref{momentlem} is proven.
\cqfd

\bigskip

\begin{remark}
Assuming $\alpha <2$, the definition of the "strange" function $q^0$ becomes:
$$ \bigl\vert \nabla \bigl( \eps^2 (1-\psi^\eps) h_2^\eps \bigr) \bigr\vert^2 Z \partial_Z (u^\eps-u^0)
\rightharpoonup  q^0 \  \mathrm{weakly\  in}\  L^2(\omega;H^{-2}(0,1)). $$ 
Indeed, if $\alpha <2$, $\lim \bigl(\eps^2 p^\eps \nabla ( \eps^2 \vert (1-\psi^\eps)h_2^\eps\vert^2 ) \bigr)
= \lim \bigl( \eps^2 p^\eps  (1-\psi^\eps)^2 \nabla ( \eps^2 \vert h_2^\eps\vert^2)  \bigr) = \int_{\mathbb{T}^{d-1}} \eta^0 \nabla_X ( \vert h_2 \vert^2) dX$
where $\eta^0$ is the two-scale limit of the sequence $(\eps^2 p^\eps  (1-\psi^\eps)^2 )$.
This limit is computable without assuming the strong $L^2$-convergence of $(\eps^2 p^\eps)$ as in \cite{BaCh}.
A first proof is in \cite{Mike}.
Here we simply use the two-scale convergence.
With the slight modification of the proof of Lemma \ref{pressurelemma}, we obtain
$$ \lim_{\eps \to 0} \Bigl( \int_\Omega \eps^2p^\eps (1-\psi^\eps) \, \div_X \phi^\eps \, dx \, dZ - \int_\Omega \eps^2p^\eps \eps^{2-\alpha} \nabla \psi^\eps \, \phi^\eps \, dx \, dZ \Bigr)  = 0,$$
for any admissible two-scale test function $\phi^\eps$.
Since $\alpha <2$, it follows that $\nabla_X \eta^0 =0$.
Thus, by the periodicity of $h_2$,
$$ \int_{\mathbb{T}^{d-1}} \eta^0 \nabla_X ( \vert h_2 \vert^2) dX = \eta^0 \int_{\mathbb{T}^{d-1}}  \nabla_X ( \vert h_2 \vert^2) dX =0 .$$
\end{remark}

%%%%%%%%%%%%%%%%%%%%%%%%%%%
%%%%%%%%%%%%%REYNOLDS
%%%%%%%%%%%%%%%%%%%%%%%%%%%

\subsection{A modified Reynolds approximation: proof of Theorem \ref{main}}

%%%%%%%%%%%%%%%%%%%%%%%%%%%

Let us check that the system formed by the equations presented in Lemmas \ref{divlem} and \ref{momentlem} corresponds to a "modified Reynolds'' system. 
The limit momentum equation is (see lemma~\ref{momentlem})
$$ N_\psi Z \partial_Z u^0 - \partial_{ZZ}^2 u^0 + h_1^2 \nabla_x p^0 + q^0 =0 . $$
We integrate it in the variable $Z$ (for each fixed $x$) and find, with the boundary conditions:
\begin{eqnarray}\label{vitesse-generalise}
&& u^0(x,Z) = \left( \int_0^Z \int_0^s e^{N_\psi(x)(s^2-t^2)/2}(\nabla_x p^0(x) + \frac{q^0(x,t)}{h_1(x)^2})\, dtds 
\right.
\nonumber \\
&& \qquad \qquad
\left.
- \int_0^1 \int_0^s e^{N_\psi(x)(s^2-t^2)/2}(\nabla_x p^0(x) + \frac{q^0(x,t)}{h_1(x)^2})\, dtds \frac{\int_0^Z e^{N_\psi(x) s^2/2}\, ds}{\int_0^1 e^{N_\psi(x) s^2/2}\, ds} \right)  h_1(x)^2 
\nonumber \\
&&\qquad \qquad
 + \left( 1-\frac{\int_0^Z e^{N_\psi(x) s^2/2}\, ds}{\int_0^1 e^{N_\psi(x) s^2/2}\, ds} \right)U_b.
\end{eqnarray}
    Next, integrating with respect to $Z$ the divergence equation (see lemma~\ref{divlem})
$$\div_x (h_1 u^0) + \dZ{(w_1+L-Z \nabla_x h_1 \cdot u^0)}=0,$$ and taking into account that the velocity $w_1+L- Z \nabla_x h_1\cdot u^0$ cancels for $Z=0$ and for $Z=1$ we obtain $$\div_x \left( \int_0^1 h_1 u^0dZ \right) = 0.$$ 
    With the previous expression \eqref{vitesse-generalise} for the velocity the following pressure equation is obtained
$$\div_x \left( h_1^3 A \nabla_x {p^0}\right) = \div_x \left(  h_1 (BU_b-Q^0) \right)$$
where $A$ and $B$ are two functions defined in $\omega$ by
\begin{equation*}
A(x) = \frac{12}{N_\psi(x)} \left( e^{N_\psi(x)/2} \int_0^1 e^{-N_\psi(x) t^2/2}\, dt - 1 \right) - \frac{12}{N_\psi(x)}(e^{N_\psi(x)/2}-1) \frac{\int_0^1 \int_0^s e^{N_\psi(x)(s^2-t^2)/2}\, dtds}{\int_0^1 e^{N_\psi(x) s^2/2}\, ds},
\end{equation*}
\begin{equation*}
B(x) = \frac{1}{N_\psi(x)}(e^{N_\psi(x)/2}-1) \frac{1}{\int_0^1 e^{N_\psi(x) s^2/2}\, ds}, \quad
x\in \omega,
\end{equation*}
assuming that function $y \mapsto (e^{y/2}-1)/y$ is continuously extended to $\mathbb{R}_+$ and $Q^0$ is defined by
\begin{eqnarray*} 
&& Q^0(x) = \int_0^1 \left( \int_0^Z \int_0^s e^{N_\psi(x)(s^2-t^2)/2}q^0(x,t)\, dtds 
\right.
\nonumber \\
&& \qquad \qquad
\left.
- \int_0^1 \int_0^s e^{N_\psi(x)(s^2-t^2)/2} q^0(x,t)\, dtds \frac{\int_0^Z e^{N_\psi(x) s^2/2}\, ds}{\int_0^1 e^{N_\psi(x) s^2/2}\, ds} \right) dZ, \quad x \in \omega.
\end{eqnarray*}
Part {\it (i)} of Theorem \ref{main} is proven.
Part {\it (ii)} is straightforward since the limit momentum equation is now given by:
$$ N_\psi Z \partial_Z u^0 - \partial_{ZZ}^2 u^0 + h_1^2 \nabla_x p^0 =0 . $$
We thus have $q^0=0$ and then $Q^0=0$.

%%%%%%%%%%%%%%%%%%%%%%%%%%%%%%%%%%%%%%%%%%%%%%%%%%%%%%%%%%%%%%%%%%%%%%%%%%
%%%%%%%%%%%%%%%%%%%%%%%%%%%%%%%%%%%%%%%%%%%%%%%%%%%%%%%%%%%%%%%%%%%%%%%%%%
%%%%%%%%%%%%%%%%%%%%%%%%%%%%%%%%%%%%%%%%%%%%%%%%%%%%%%%%%%%%%%%%%%%%%%%%%%
\section{Some numerical illustrations}\label{SectNum}
%%%%%%%%%%%%%%%%%%%%%%%%%%%%%%%%%%%%%%%%%%%%%%%%%%%%%%%%%%%%%%%%%%%%%%%%%%
%%%%%%%%%%%%%%%%%%%%%%%%%%%%%%%%%%%%%%%%%%%%%%%%%%%%%%%%%%%%%%%%%%%%%%%%%%
%%%%%%%%%%%%%%%%%%%%%%%%%%%%%%%%%%%%%%%%%%%%%%%%%%%%%%%%%%%%%%%%%%%%%%%%%%

In this section, we present some numerical results.
The goal is to quantify the effect of roughness on a portion of a surface.
So we compare the solutions of the classical Reynolds equation (corresponding to the smooth domain) and those obtained from the equation with partial roughness, see Theorem~\ref{main} on page~\pageref{main} and Remark \ref{base-num} after.
\par
\noindent
Physically, we are in the framework of lubrication.
Consequently, we are interestd in the pressure forces in a flow between two surfaces whose relative velocity is imposed.
As usual, all the geometry of the gap between the two surfaces is taking into account by the upper surface (it corresponds to the height~$h_1$) and the relative velocity is taking account by the lower surface (its imposed velocity is denoted by~$U_b$, while the upper surface is at rest).
Moreover, we know that the problem is well posed if and only if we give boundary conditions: either by imposing Dirichlet condition on the pressure, either by imposing Neumann type flow.
A physically interesting situation is the following: we give a Neumann condition at the ``entrance'' of the channel (that is a flux denoted by~$Q_e$), and we impose a pressure at the other boundaries.
\par
\noindent
More precisly, for the next simulations we use the following data (see also Fig.~\ref{3dflow}):
\begin{equation*}
\omega=[0,1]\times[0,1],
\quad
h_1(x,y)=(2x-1)^2+0.5,
\quad
U_b=(1,0)
\quad \text{and} \quad
Q_e=0.5
\end{equation*}

%%%%%%%%%%%%%%%%%%%%%%%%%%%%%%%%%%%%%%%%%%%%%%%%%%%%%%%%%%%%%%%%%%%%%%%%%%
\subsection{Without rugosities}
%%%%%%%%%%%%%%%%%%%%%%%%%%%%%%%%%%%%%%%%%%%%%%%%%%%%%%%%%%%%%%%%%%%%%%%%%%

\begin{figure}[htbp]
\centering
{\psfrag{Q0}{Flux}\psfrag{Ub}{Bottom velocity}\psfrag{H}{Height}\psfrag{Hbis}{(with possible rugosities)}
\includegraphics[height=4cm]{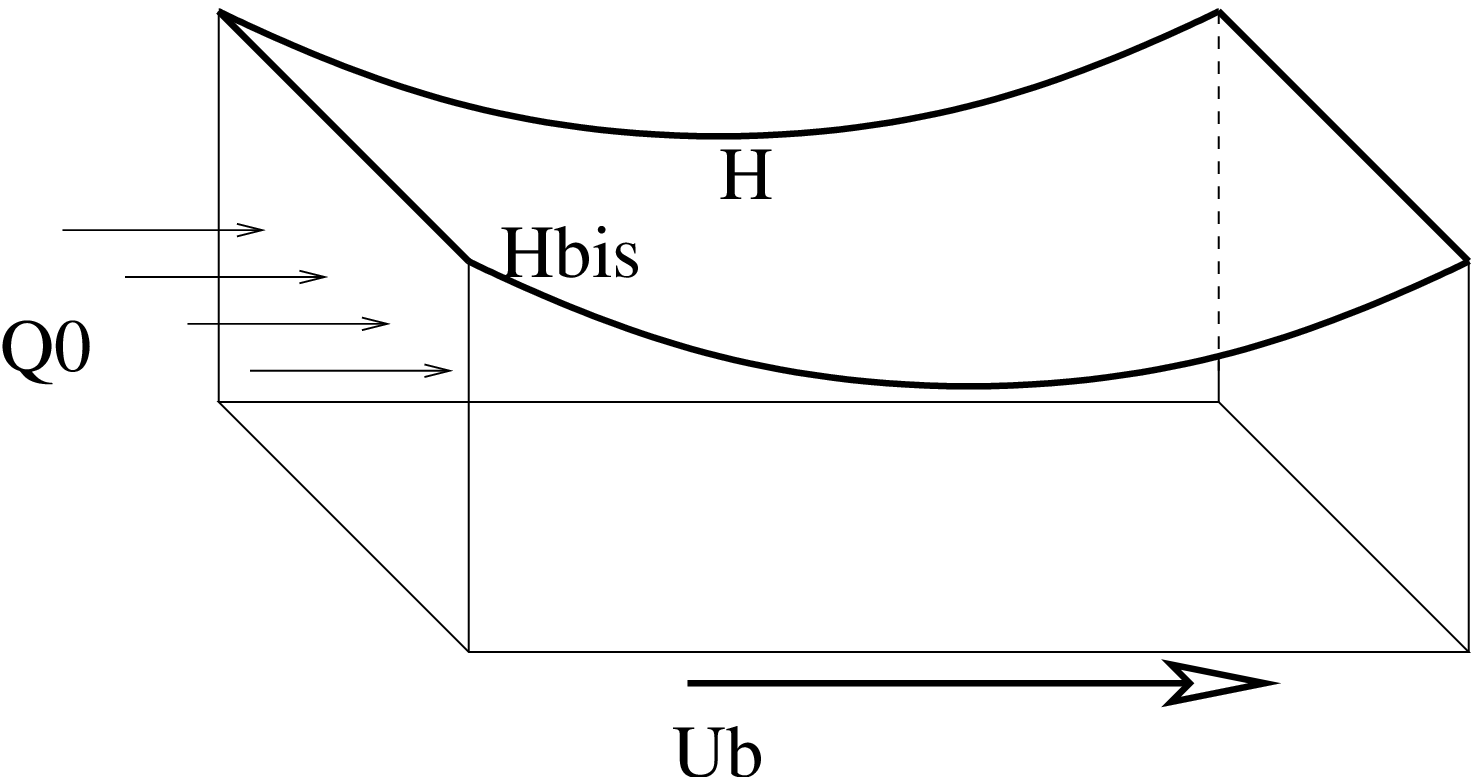}
}
\includegraphics[height=6cm]{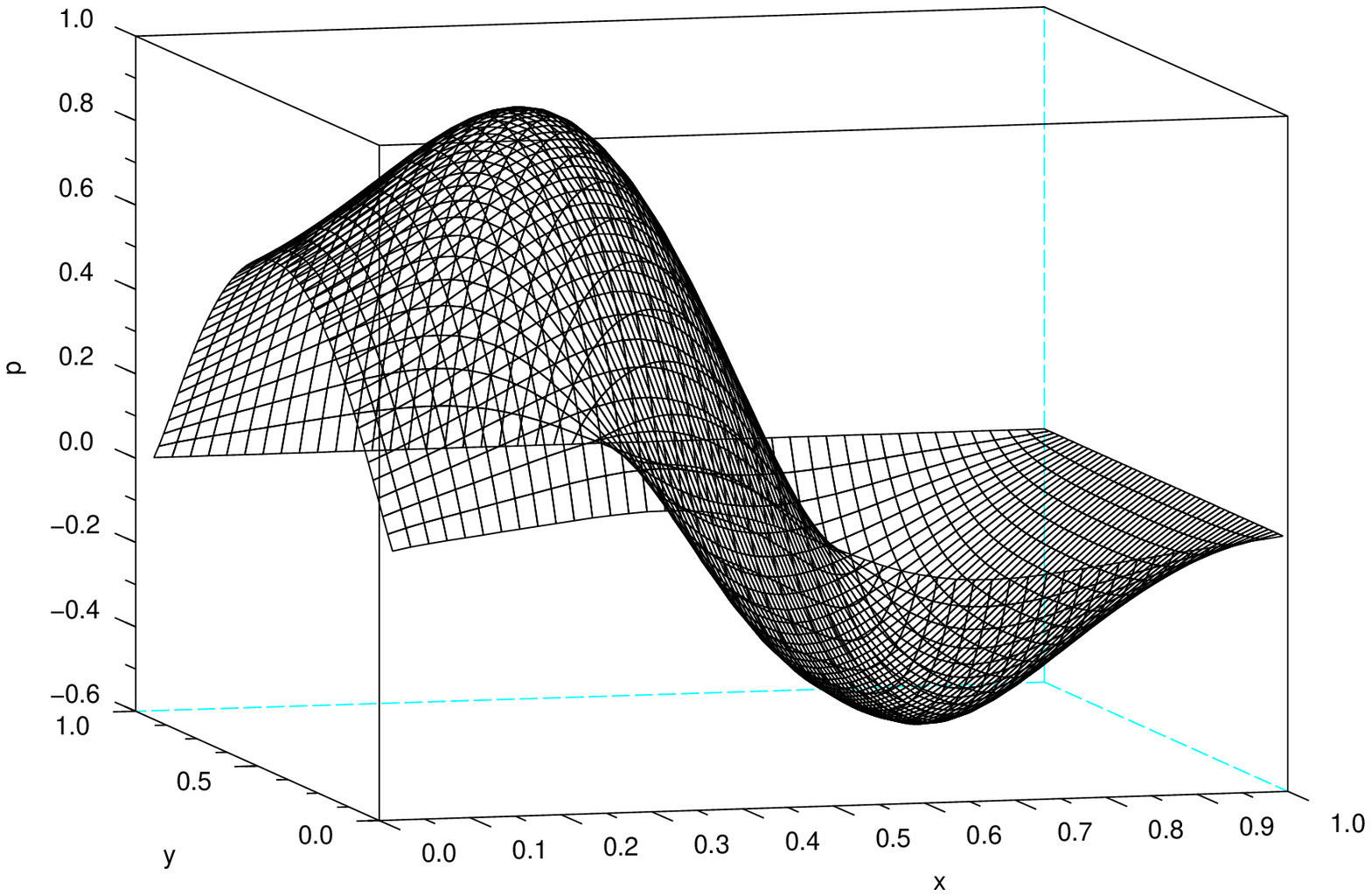}
\caption{On the left figure, we draw the three-dimensional domain and indicate the main data.
On the right figure, we plot the pressure profile without rugosity, that is the solution to the classical Reynolds equation.}\label{3dflow}
\end{figure}

To situate the next results, the first figure that we present (see Fig.~\ref{3dflow}, on right) corresponds to the case of flow in a smooth domain.
In that case, we solve the classical Reynolds equation, see the Remark~\ref{rem-1202}.
We use the Freefem++ language (see http://www.freefem.org/ff++/ ) and consequently a variational form of the Reynolds equation: for all test function $q$ vanishing on the boudaries $\{x=1\}$, $\{y=0\}$ and $\{y=1\}$, we get
\begin{equation*}
  \int_\omega \frac{h_1^3}{12 } \nabla p \cdot \nabla q - \int_\omega \frac{h_1}{2} U_b\cdot \nabla q + \int_{\{x=0\}}Q_e q = 0.
\end{equation*}
Finally, we use the Scilab software (http://www.scilab.org/ ) to view results.

%%%%%%%%%%%%%%%%%%%%%%%%%%%%%%%%%%%%%%%%%%%%%%%%%%%%%%%%%%%%%%%%%%%%%%%%%%
\subsection{Some rugosities}
%%%%%%%%%%%%%%%%%%%%%%%%%%%%%%%%%%%%%%%%%%%%%%%%%%%%%%%%%%%%%%%%%%%%%%%%%%

\begin{figure}[htbp]
\centering
{\psfrag{R}{Rugosity part}\psfrag{L}{Smooth part}
\includegraphics[height=3.5cm]{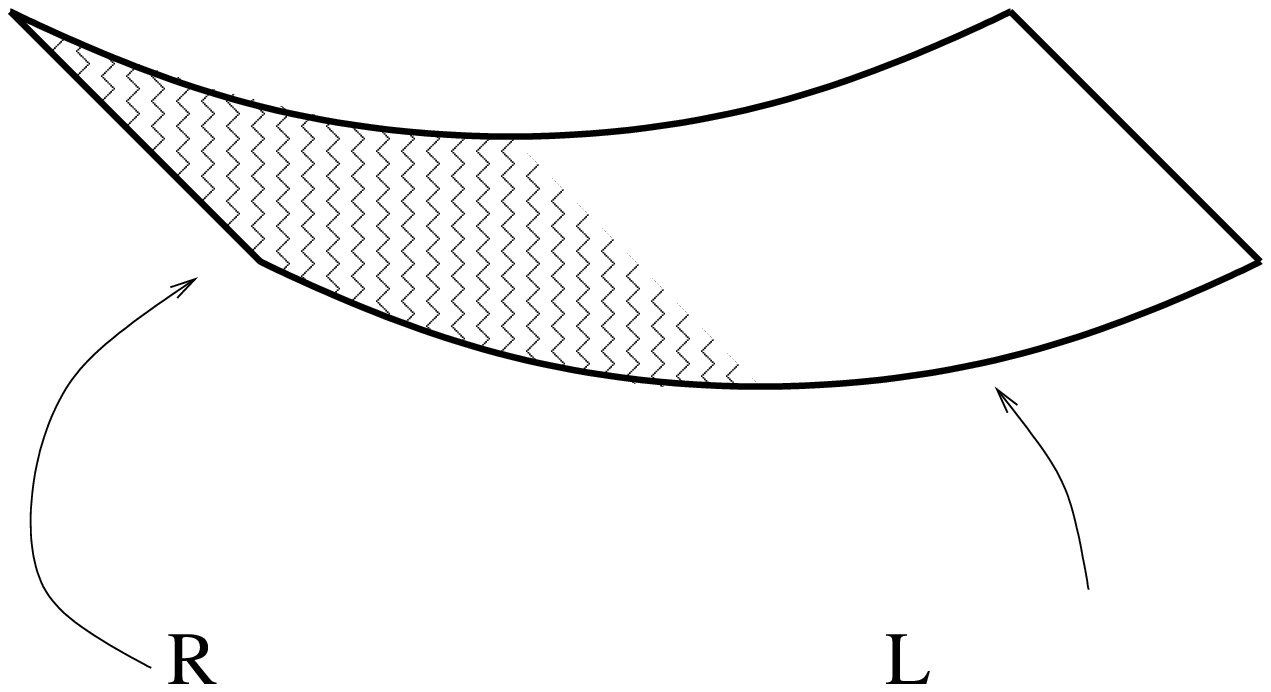}
}
\hspace{1cm}
\includegraphics[height=6cm]{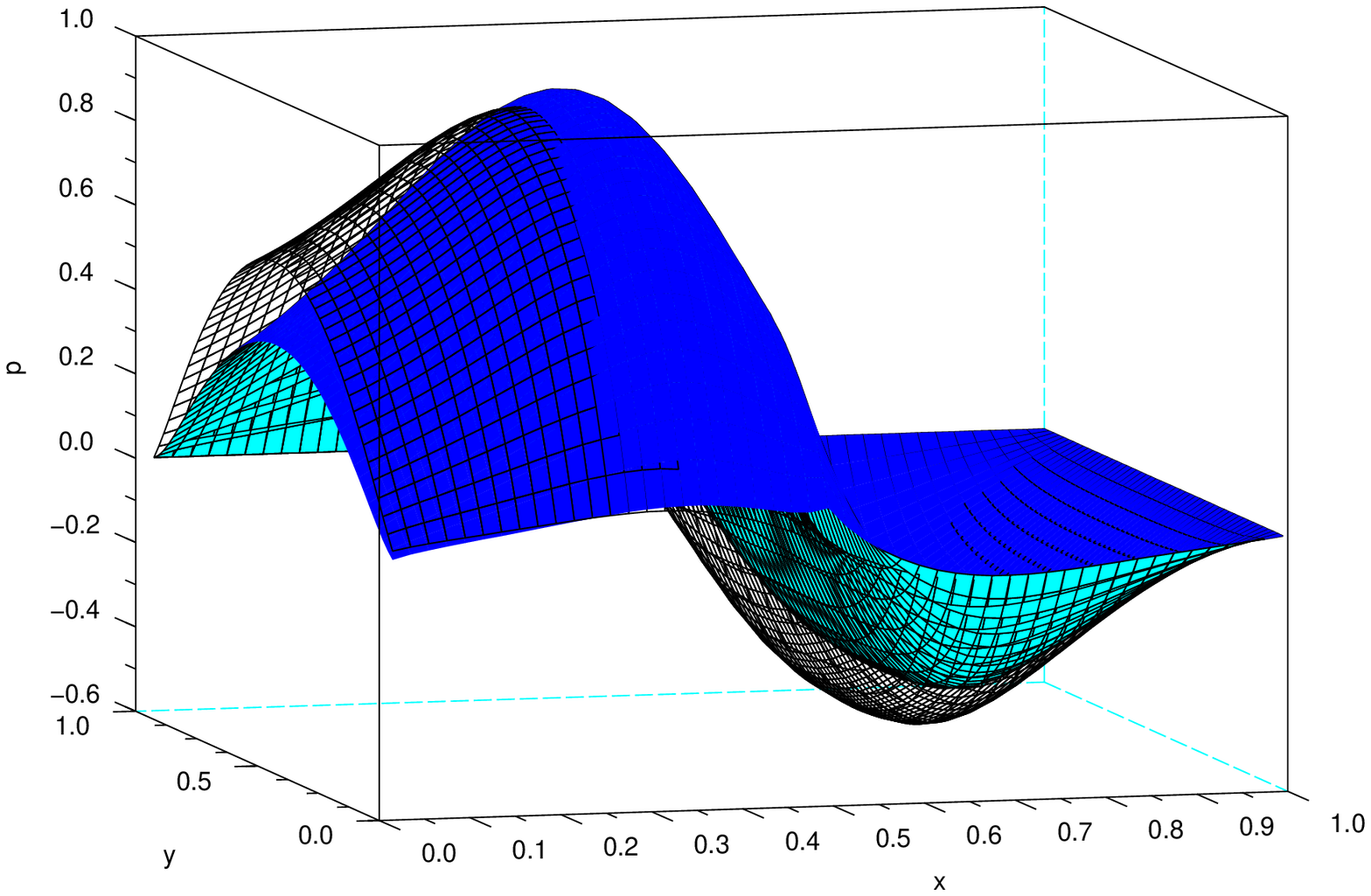}
%On a aussi a disposition une coupe : fig-rug1
\caption{On the left figure, we indicate the part of the top boundary which is rough.
On the right figure, we compare the pressure profile without rugosity (grid) and the pressure profile obtain with the corresponding rugosities (plain).% To better highlight the effects, we only present a half domain.
}\label{3dflow-rug1}
\end{figure}

\begin{figure}[htbp]
\centering
{\psfrag{L}{Rugosity part}\psfrag{R}{Smooth part}
\includegraphics[height=3.5cm]{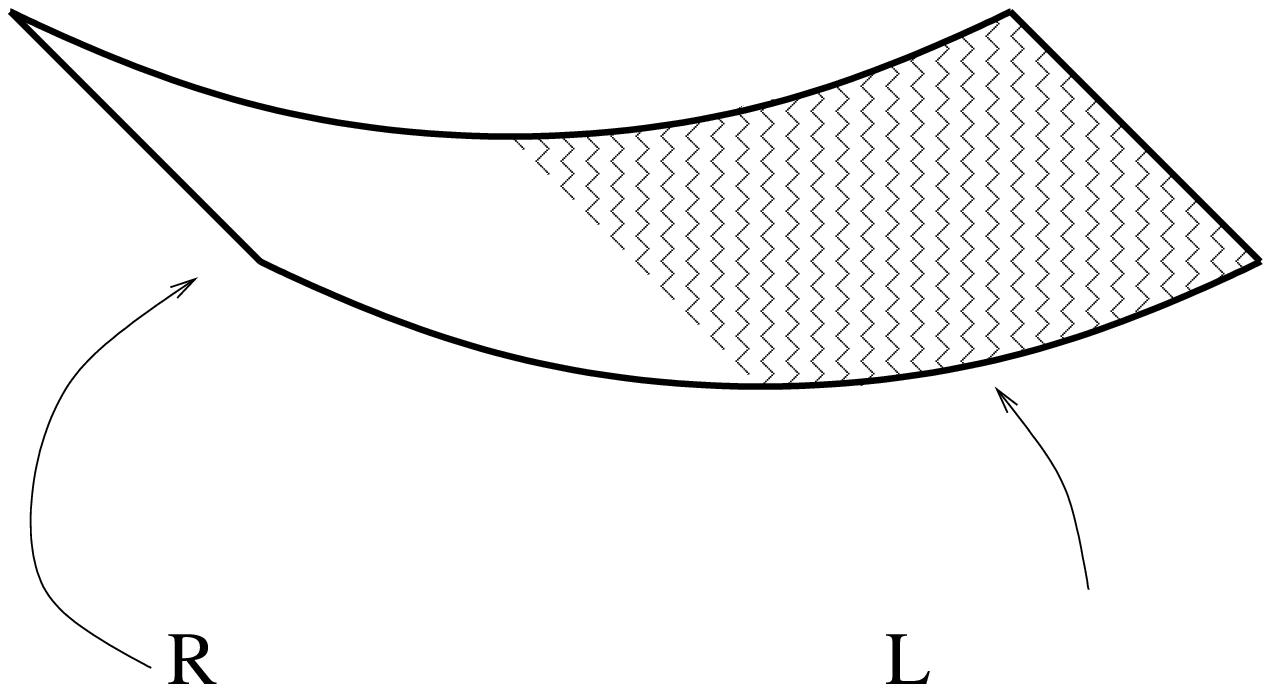}
}
\hspace{1cm}
\includegraphics[height=6cm]{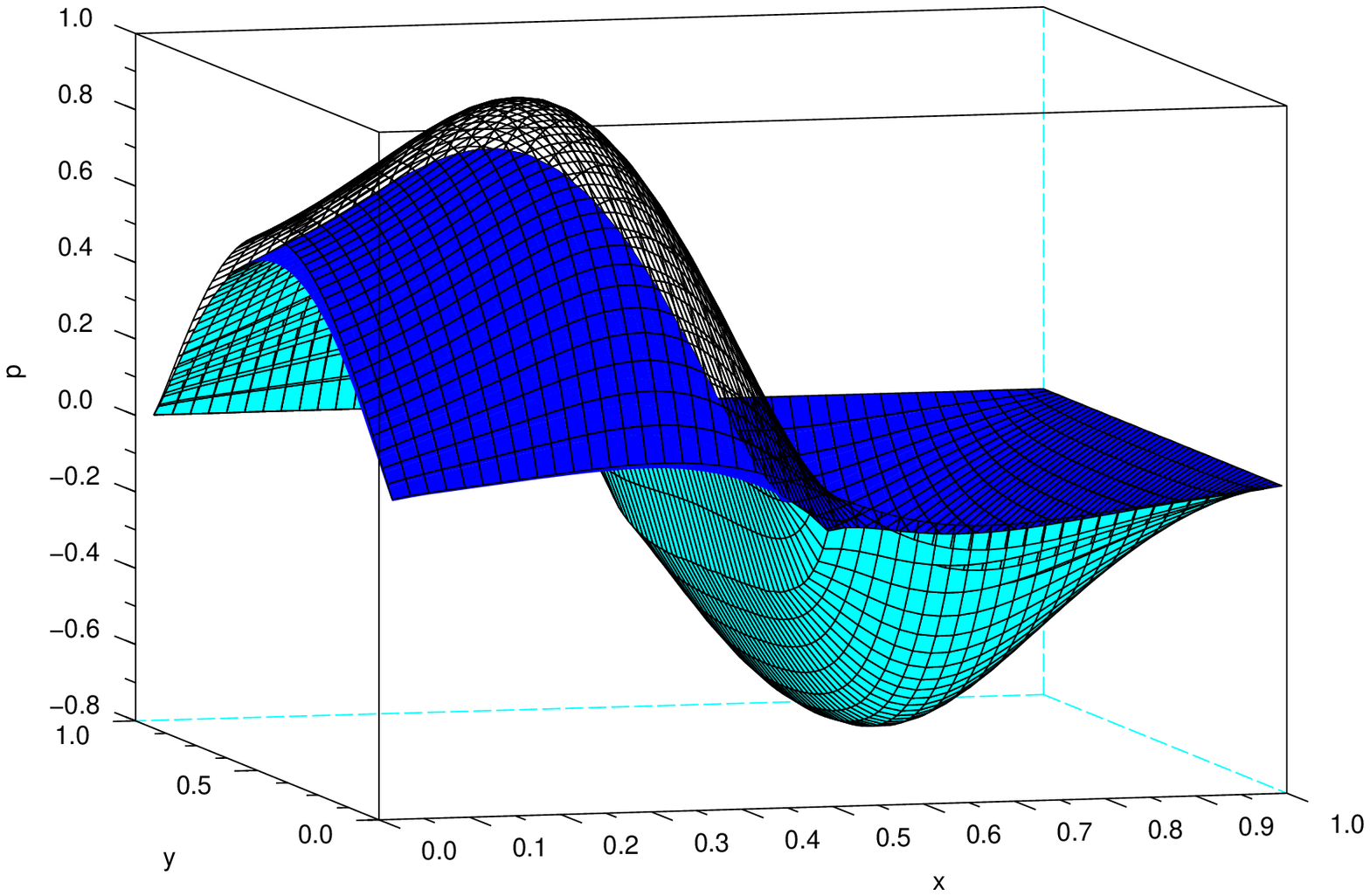}
\caption{On the left figure, we indicate the part of the top boundary which is rough.
On the right figure, we compare the pressure profile without rugosity (grid) and the pressure profile obtain with the corresponding rugosities (plain).}\label{3dflow-rug2}
\end{figure}

\begin{figure}[htbp]
\centering
{\psfrag{L}{Rugosity part}\psfrag{R}{Smooth part}
\includegraphics[height=3.5cm]{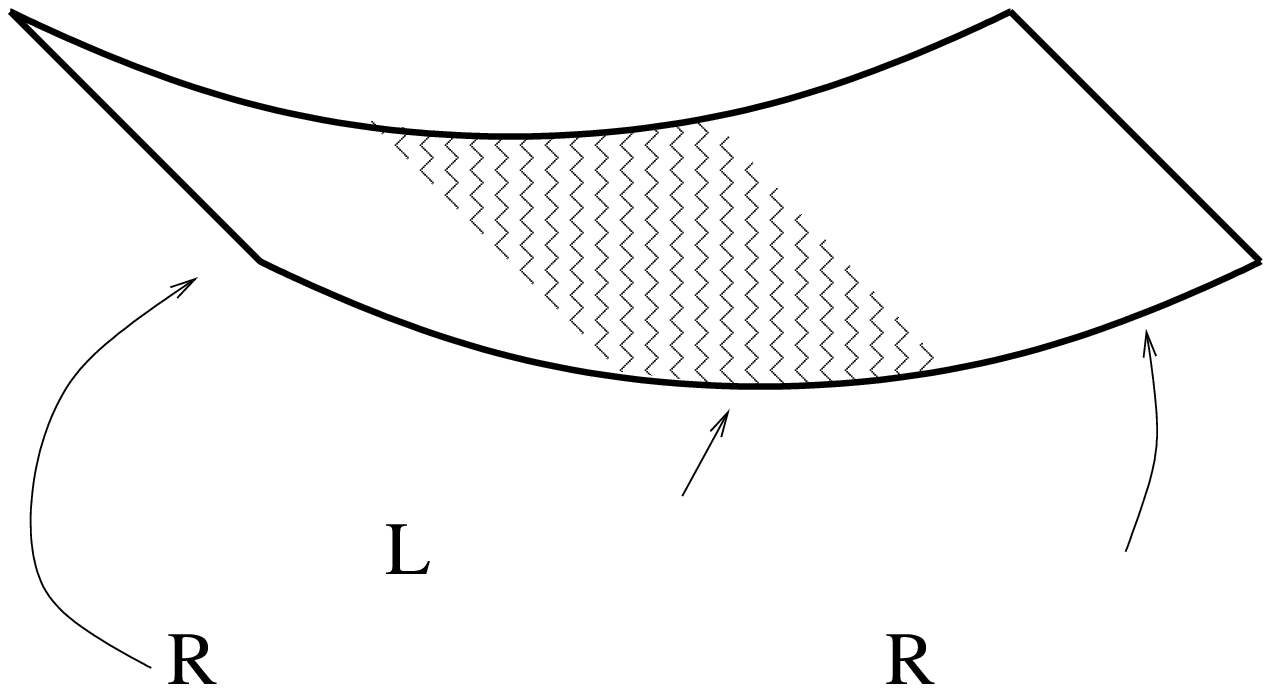}
}
\hspace{1cm}
\includegraphics[height=6cm]{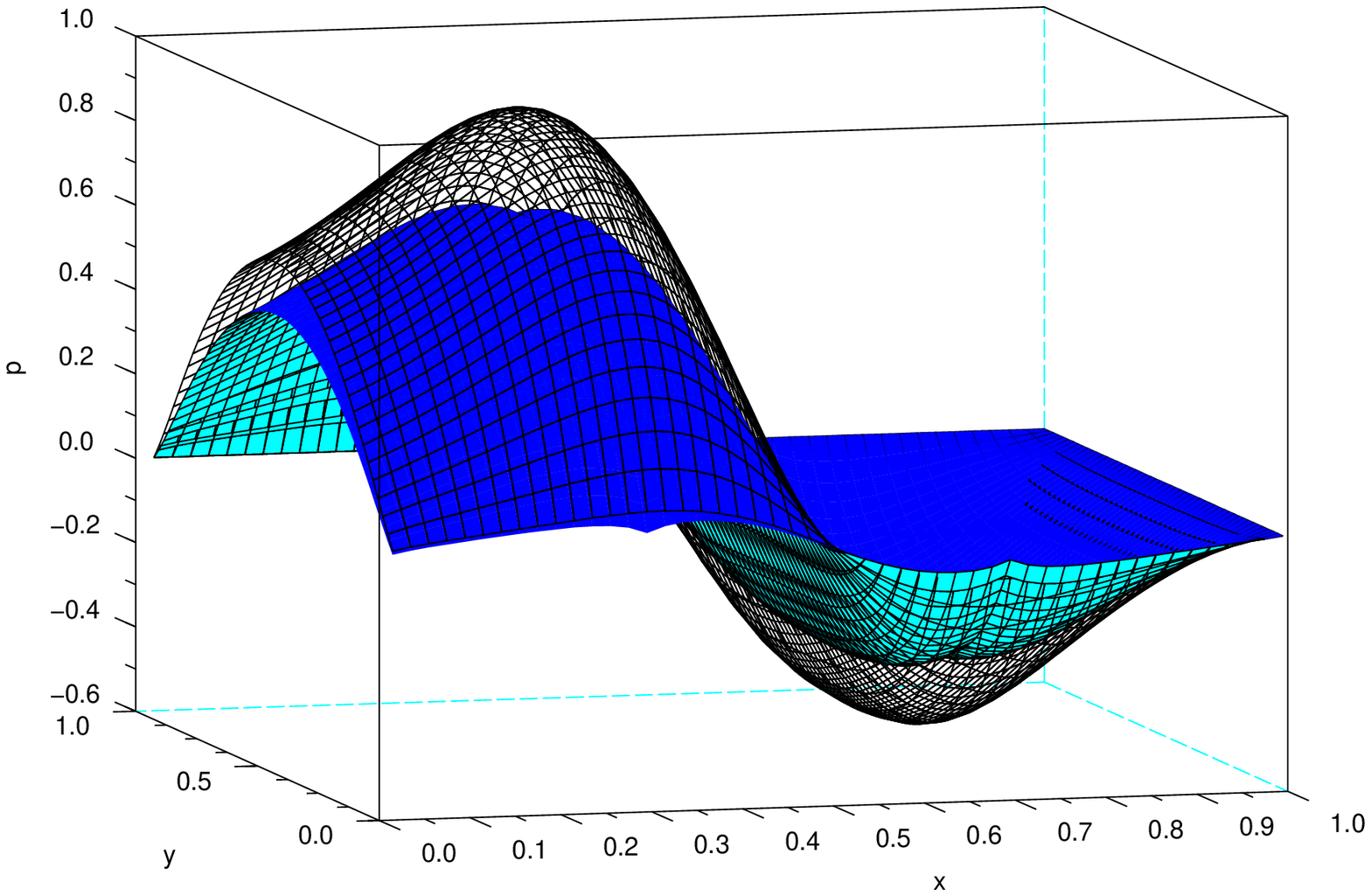}
\caption{On the left figure, we indicate the part of the top boundary which is rough.
On the right figure, we compare the pressure profile without rugosity (grid) and the pressure profile obtain with the corresponding rugosities (plain).}\label{3dflow-rug3}
\end{figure}

From the results of Theorem~\ref{main}, it seems simple to simulate the rugosities effets on the upper surface.
Morally, it is sufficient to specify first: on which area of the surface roughness ago, and secondly: quantify the roughness using their form.
These data allow to evaluate the function~$N_\psi$, and therefore the coefficients~$A$ and~$B$ involved in the statement of the Theorem~\ref{main}.
For standard rugosities, for instance described by a function of the form $\alpha \cos(\beta X)$, it is reasonable to take the following examplarity values 
\begin{equation*}
\begin{aligned}
& N_\psi(x,y) = 0, \quad A(x,y) = 1 \quad \text{and} \quad B=0.5 \quad \text{on smooth part,}\\
& N_\psi(x,y) = 2, \quad A(x,y) = 1.08696 \quad \text{and} \quad B=0.58739 \quad \text{on rough part.}\\
\end{aligned}
\end{equation*}
We then present three cases where the domain~$\omega$ is partially rough (see Fig.~\ref{3dflow-rug1}, \ref{3dflow-rug2} and~\ref{3dflow-rug3}). 
We can note two remarkable points.
First, the effect induced by the roughness is not a local effect.
Thus, although roughness occupy only a portion of the domain, then the pressure is disrupted throughout the domain, including upstream roughness.
Secondly, the effect of roughness is quite significant even if these wrinkles are present only a small part of the domain (see Fig.~\ref{3dflow-rug3}).

\noindent {\bf Acknowledgments -}
The second author of this work has been partially supported by the ANR project n° ANR-08-JCJC-0104-01 : RUGO (Analyse et calcul des effets de rugosités sur les écoulements).

%%%%%%%%%%%%%%%%%%%%%%%%%%%%%%%%%%%%%%%%%%%%%%%%%%%%%%%%%%%%%%%%%%%%%%%%%%

%%%%%%%%%%%%%%%%%%%%%%%%%%%%%%%%%%%%%%%%%%%%%%%%%%%%%%%%%%%%%%%%%%%%%%%%%%


\begin{thebibliography}{00}

%%%%%%%%%%%%%%%%%%%%%%%%%%%%%%%%%%%%%%%%%%%%%%%%%%%%%%%%%%%%%%%%%%%%%%%%%%


\bibitem{AcPiVa}
 \sc Y. Achdou, O. Pironneau, F. Valentin.
 \rm Effective boundary conditions for laminar flows over periodic rough boundaries.
 \it J. Comput. Phys., \rm {147-1}, 187--218, (1998).

 \bibitem{A}
\sc G.~Allaire.
\rm Homogenization and two-scale convergence.
{\em SIAM J. Math. Anal.}, 23:1482--1518, (1992).

 \bibitem{BaCh1}
 \sc G. Bayada, M. Chambat.
 \rm New models in the theory of the hydrodynamic lubrication of rough
 surfaces.
 \it Trans. of the AMS J. of Trib.,  \rm {\bf 110}, 402-407, (1988).
 
  \bibitem{BaCh}
 \sc G. Bayada, M. Chambat.
 \rm Homogenization of the Stokes system in a thin film flow with rapidly varying thickness.
 \it  RAIRO Model. Math. Anal. Numer, 
 \rm {\bf 23}, no. 2, 205--234, (1989).
 

 \bibitem{BeChCi}
 \sc N. Benhaboucha, M. Chambat, I. Ciuperca.
 \rm Asymptotic behaviour of pressure and stresses in a thin film flow
 with a rough boundary.
 \it Quart. of Appl. Math. \rm no. 2, 369--400, (2006).
 
  \bibitem{BrChChCoGi}
 \sc D. Bresch, C. Choquet, L. Chupin, T. Colin,  M. Gisclon.
 \rm Roughness-Induced Effect at Main order on the Reynolds Approximation.
 {\it SIAM Multiscale Model. Simul.}, to appear (2010).
  
  \bibitem{BhuBla91}
 \sc B.~Bushan and G.~Blackman.
\newblock \rm Atomic force microscopy of magnetic rigid disks disks and sliders and
  its applications to tribology.
\newblock {\em ASME J. Tribol.}, 113:452--457, (1991).


 \bibitem{CiMu} 
\sc D. Cioranescu, F. Murat. 
\rm Un terme \'etrange venu d'ailleurs.
\rm [A strange term brought from somewhere else].
\it Nonlinear partial differential equations and their applications.
\rm Coll\`ege de France Seminar, Vol. II (Paris, 1979/1980),
pp. 98--138, 389--390, Res. Notes in Math., 60, Pitman, Boston, Mass.-London, (1982).

 \bibitem{Dri}
 \sc H. Dridi.
 \rm Comportement asymptotique des \'equations de Navier-Stokes dans des domaines
applatis. 
 {\it Bull. Sc. Math.}, 106, 369--385 (1982). 
 
 \bibitem{Dy}
 \sc A. Dyson.
 \rm Hydrodynamic lubrication of rough surface - a review work. 
 \it Proceedings of the 4th Leeds-Lyon Symposium on surfaces roughness on lubrication,
  \rm 61--69, (1977).

\bibitem{El}
 \sc H.-G. Elrod
 \rm A review of theories for the fluid dynamic effects of roughness on laminar lubricating films. 
 \it Proceedings of the 4th Leeds-Lyon Symposium on surfaces roughness on lubrication,
  \rm 11--26, (1977).

  \bibitem{ChaFit}
 \sc T. Fitzgerald, R. Farrell, N. Petkov, C. Bolger, M. Shaw, J. Charpin, J.P. Gleeson, J. Holmes, and M. Morris. 
\newblock  \rm A Study on the Combined Effects of Solvent Evaporation and Polymer Flow upon Block Copolymer Self-Assembly and Alignement on Topographic Patterns.
\it Langmuir,
\rm 25 (23), 13551--13560

\bibitem{HeZhu97}
 \sc  L.~He and J.~Zhu.
\newblock \rm The fractal character of processed metal surfaces.
\newblock {\em Wear}, 208:17--24, (1997).


\bibitem{JaMi2}
 \sc W. J\"ager, A. Mikeli\'c.
 \rm Couette flows over a rough boundary and drag reduction. 
 \it Comm. Math. Phys, 
 \rm (232-3), 429--455, (2003).
 
 \bibitem{Kou}
\sc M.-M Koura and M.-A. Omar.
\rm The effects of  surface parameters on friction.
\it Wear, 
\rm 73(2): 235--246, (1981).

\bibitem{Mike}
 \sc A. Mikeli\'c.
 \rm Remark on the result on homogenization in hydrodynamical lubrication by G. Bayada and M. Chambat.
 \it RAIRO Model. Math. Anal. Numer, 
 \rm (25-3), 363--370, (1991).

\bibitem{My}
\sc N.-O. Myers.
\rm Characterization of surface roughness.
\it Wear, 
\rm 5(3): 182--189, (1962).

\bibitem{Ng}
\sc G. Nguetseng.
\rm A general convergence result for a functional related to the theory
  of homogenization.
{\em SIAM J. Math. Anal.}, 20:608--623, (1989).

\bibitem{PPL}
\sc F. Plourabou\'e, M. Prat, N. Letalleur.
\rm Sliding lubricated anisotropic rough surfaces. 
{\it Physical Review E}, 64, 011202, (2001).
 
 \bibitem{Rey86}
 \sc O.~Reynolds.
\newblock \rm On the theory of lubrication and its application to mr beauchamp
  tower's experiment, including an experimental determination of the viscosity
  of olive oil.
\newblock {\em Proc. Roy. Soc. London}, 40:191--203, (1886).



\bibitem{Wiel}
\sc W. Wieleba.
\rm The statistical correlation of the coefficient of friction and wear rate of PTFE composites with steel counterface roughness and hardness.
\it  Wear, 
\rm 252(9-10): 719--729, (2002).

\bibitem{Zhik04}
 \sc  V.-V.~Zhikov.
\newblock \rm On two-scale convergence.
\newblock {\em J. Math. Sci.}, 120(3):1328--1352, (2004).


\bibitem{ZhoLeu93}
 \sc  G.~Zhou, M.-C. Leu, and D.~Blackmore.
\newblock \rm Fractal geometry model for wear prediction.
\newblock {\em Wear}, 170:1--14, (1993).

%%%%%%%%%%%%%%%%%%%%%%%%%%%%%%%%%%%%%%%%%%%%%%%%%%%%%%%%%%%%%%%%%%%%%%%%%%

\end{thebibliography}
\end{document}